\documentclass[12pt]{amsart}

\pdfoutput=1

\usepackage[a4paper]{geometry}
\usepackage{amsmath,amssymb,amsthm,amsfonts,color}
\usepackage{comment}
\usepackage{hyperref}
\usepackage{stmaryrd}
\usepackage{tikz}
\usepackage{bm}

\usepackage[shortlabels]{enumitem}
\usepackage{refcount}

\usepackage{threeparttable}
\usepackage{caption}
\theoremstyle{definition}
\newtheorem{Def}{Definition}

\newtheorem{Thm}[Def]{Theorem}
\newtheorem{Lem}[Def]{Lemma}

\newtheorem{Rem}[Def]{Remark}
\newtheorem{Prop}[Def]{Proposition}

\newtheorem{Cor}[Def]{Corollary}
\newtheorem{Claim}{Claim}
\newtheorem{Subclaim}{Subclaim}
\newtheorem{Q}{Question}

\newtheorem{Thm?}[Def]{Theorem?}

\newcommand{\ZFC}{\ensuremath{\mathsf{ZFC}}}

\newcommand{\Qp}[1]{\left\llbracket #1 \right\rrbracket}

\newcommand{\M}{\text{M}}
\newcommand{\Lp}{\text{Lp}}

\renewcommand{\L}{\text{L}}

\begin{document}


%\primarydata{03-06, 03C95, 03E40}
%\secondarydata{03E57}


\title[Boolean-valued second-order logic revisited]{Boolean-valued second-order logic revisited}
%\title{Universally Baire sets in $2^{\kappa}$}
%\author{Daisuke Ikegami and Matteo Viale}
\date{\today}

\author[D.\ Ikegami]{Daisuke Ikegami}
\address[D.\ Ikegami]{Institute of Logic and Cognition, Department of Philosophy, Sun Yat-sen University, Xichang Hall 602, 135 Xingang west street, Guangzhou, 510275 CHINA}

\email[D.\ Ikegami]{\href{mailto:ikegami@mail.sysu.edu.cn}{ikegami@mail.sysu.edu.cn}}

\thanks{The author would like to thank Jouko V\"{a}\"{a}n\"{a}nen, Toshimichi Usuba, and Nam Trang for discussions and helpful comments on the topic of this paper.}

\begin{abstract}
Following the paper~\cite{BVSOL} by V\"{a}\"{a}n\"{a}nen and the author, we continue to investigate on the difference between Boolean-valued second-order logic and full second-order logic. 
We show that the compactness number of Boolean-valued second-order logic is equal to $\omega_1$ if there are proper class many Woodin cardinals. This contrasts the result by Magidor~\cite{compactness2ndordermagidor} that the compactness number of full second-order logic is the least extendible cardinal. We also introduce the inner model $C^{2b}$ constructed from Boolean-valued second-order logic using the construction of G\"{o}del's Constructible Universe L. We show that $C^{2b}$ is the least inner model of $\mathsf{ZFC}$ closed under $\M_n^{\#}$ operators for all $n < \omega$, and that $C^{2b}$ enjoys various nice properties as G\"{o}del's L does, assuming that Projective Determinacy holds in any set generic extension. This contrasts the result by Myhill and Scott~\cite{Myhill-Scott} that the inner model constructed from full second-order logic is equal to HOD, the class of all hereditarily ordinal definable sets.
\end{abstract}

\maketitle

\section{Introduction}

%Full second-order logic (i.e., second-order logic with full structures where the second-order variables range over all subsets and relations of the first-order domain of the structures

In second-order logic, there are two major ways to interpret second-order quantifiers. 
In so-called full second-order logic, the second-order variables range over all subsets and relations of the first-order domain of the structures. 
In so-called Henkin second-order logic, the second-order variables range over a fixed set of subsets and relations of the first-order domain which satisfies Comprehension schema. 
Full second-order logic is very complicated and difficult to analyze while Henkin second-order logic is as simple as first-order logic and it enjoys nice properties such as completeness theorem and compactness theorem. 

In \cite{BVSOL}, V\"{a}\"{a}n\"{a}nen and the author introduced Boolean-valued second-order logic where the second-order variables range over all Boolean-valued subsets and relations of the first-order domain of the structures. Under the existence of large cardinals and some assumptions on Woodin's $\Omega$-logic, we showed that Boolean-valued second-order logic is simpler and easier to treat than full second-order logic in various aspects. %the aspects of validity, completeness, compactness numbers, Hanf numbers, and L\"{o}wenheim-Skolem numbers. 

In this paper, we summarize the known results and present new results on the difference between Boolean-valued second-order logic and full second-order logic with respect to the following notions and numbers:
\begin{itemize}
\item The complexity of the set of G\"{o}del numbers of the sentences in the logic which are valid.

\item The compactness number (i.e., the least infinite cardinal $\kappa$ such that any theory in the logic, all subsets of size less than $\kappa$ of which have a model, has a model).

\item The Hanf number (i.e., the least infinite cardinal $\kappa$ such that any sentence in the logic which has a model of cardinality at least $\kappa$ has arbitrarily large models).

\item The L\"{o}wenheim number (i.e., the least infinite cardinal $\kappa$ such that any sentence in the logic which has a model, has a model of cardinality at most $\kappa$). 

\item The inner model constructed from the logic using the construction of G\"{o}del's Constructible Universe L. 
\end{itemize}

Regarding validity, let $0^2$ be the set of  G\"{o}del numbers of the formulas which are valid in full second-order logic and $0^{2b}$ be the corresponding set in Boolean-valued second-order logic. It is known that $0^2$ is not invariant under set forcings (e.g., the Continuum Hypothesis is equivalent to the validity of a sentence in full second-order logic), and that $0^2$ is $\Pi_2$-complete in $\mathsf{ZFC}$ (\cite[Theorem~1]{full-SOL}). 
V\"{a}\"{a}n\"{a}nen and the author showed that under the existence of proper class many Woodin cardinals, the set $0^{2b}$ is invariant under set forcings and $0^{2b}$ is Turing equivalent to $0^{\Omega}$, the set of G\"{o}del number of sentences which are $\Omega$-valid in Woodin's $\Omega$-logic (\cite[Theorem~2.7, Theorem~5.1]{BVSOL}). If there are proper class many Woodin cardinals and Woodin's $\Omega$-Conjecture holds, then $0^{\Omega}$ is $\Delta_2$-definable in the language of set theory. Hence it is consistent that $0^{2b}$ is strictly simpler than $0^{\Omega}$ (\cite[Corollary~6.8]{BVSOL}). 
Furthermore, if we additionally assume that $\mathsf{AD}^+$-Conjecture holds, then $0^{\Omega}$ is definable over the structure $(H_{\mathfrak{c}^+} , \in)$ where $\mathfrak{c} = 2^{\omega}$, in other words, $0^{\Omega}$ is definable over the third-order arithmetic $\left(\omega , \wp (\omega) , \wp \bigl( \wp(\omega) \bigr) , \in \right)$ (\cite[Subsection~8.3]{Omega-logic}). Therefore, it is consistent that $0^{2b}$ is definable over the third-order arithmetic. This is optimal in the sense that $0^{\Omega}$ (or $0^{2b}$) cannot be definable over the second-order arithmetic $\bigl( \omega , \wp (\omega) , \in \bigr)$ (in fact, $0^{\Omega}$ is even more complex than the theory of the structure $\bigl(\mathrm{L}(\mathbb{R}), \in\bigr)$).

Regarding compactness, Magidor proved that the compactness number of full second-order logic is the least extendible cardinal if it exists (\cite[Theorem~4]{compactness2ndordermagidor}). 
V\"{a}\"{a}n\"{a}nen and the author proved that the compactness number of Boolean-valued second-order logic is at most the least supercompact cardinal if it exists under the existence of proper class many Woodin cardinals and the assumption that the strong $\Omega$-Conjecture holds (\cite[Theorem~9.10]{BVSOL}). 

In this paper, we show that the compactness number of Boolean-valued second-order logic is equal to $\omega_1$ under the existence of proper class many Woodin cardinals (Corollary~\ref{cor:compactness-BVSOL}).  
%we introduce generically extendible cardinals and weakly generically extendible cardinals, and show that the compactness number of Boolean-valued second-order logic is the least weakly generically extendible cardinal if it is exists, and that both the least generically extendible cardinal and the least weakly generically extendible cardinal are equal to $\omega_1$ under the existence of proper class many Wooden cardinals (Theorem~\ref{thm:L2bcompactness-number} and Corollary~\ref{cor:compactness-BVSOL}). Hence it is consistent that the compactness number of Boolean-valued second-order logic is much smaller than that of full second-order logic. 
%Furthermore, Usuba showed that the existence of generically extendible cardinals does not imply $0^{\#}$ (\cite[Theorem~1.4]{Usuba}), and that $\mathsf{ZFC}$ plus the existence of generically extendible cardinals is equiconsistent with $\mathsf{ZFC}$ plus the existence of weakly generically extendible cardinals (Corollary~\ref{cor:equiconsistency}). Hence the existence of the compactness number of Boolean-valued second-order logic does not imply $0^{\#}$ either while the existence of the compactness number of full second-order logic implies that of extendible cardinals, which is much stronger than that of $0^{\#}$. 

Regarding the Hanf numbers, Magidor proved that the Hanf number $h^2$ of full second-order logic is between the least supercompact cardinal and the least extendible cardinal if it exists (\cite{compactness2ndordermagidor}).  V\"{a}\"{a}n\"{a}nen and the author proved that the Hanf number $h^{2b}$ of Boolean-valued second-order logic is less than the least supercompact cardinal if there are proper class many Woodin cardinals and the Strong $\Omega$-Conjecture holds (\cite[Theorem~9.13]{BVSOL}). 
Recently, Usuba proved that $h^{2b}$ is equal to $\omega$ if there are proper class many Woodin cardinals (\cite{Usuba}).  %(\cite[Corollary~7.5]{Usuba}). %$h^{2b}$ is less than the least generically extendible cardinal (\cite[Proposition~7.4]{Usuba}. Since the least generically extendible cardinal is $\omega_1$ if there are proper class many Woodin cardinals, $h^{2b}$ is equal to $\omega$ if there are proper class many Woodin cardinals (\cite[Corollary~7.5]{Usuba}). Hence it is consistent that $h^{2b}$ is $\omega$ while $h^2$ is greater than the first supercompact cardinal.

Regarding the L\"{o}wenheim number, it is known that the L\"{o}wenheim number $\ell^2$ of full second-order logic is the supremum of $\Delta_2$-definable ordinals in $\mathsf{ZFC}$ (see e.g., \cite{Lowenheim-number}). V\"{a}\"{a}n\"{a}nen and the author showed that it is consistent that the L\"{o}wenheim number $\ell^{2b}$ of Boolean-valued second-order logic is less than the least Woodin cardinal (\cite[Theorem~9.17]{BVSOL}). Hence it is consistent that $\ell^{2b} < \ell^2$. 

Regarding the inner model construction, G\"{o}del constructed the Constructible Universe L by repeating to add first-order definable subsets of the structure in question. G\"{o}del's L is the smallest inner model of $\mathsf{ZFC}$. It is defined absolutely among all transitive proper class models of $\mathsf{ZF}$. G\"{o}del's L admits a fine structure and it satisfies various combinatorial principles such as $\mathsf{GCH}$, $\diamondsuit_{\kappa}$ for all regular uncountable cardinals $\kappa$, $\square_{\kappa}$ for all uncountable cardinals $\kappa$. 
If one uses a logic $\L^*$ extending first-order logic to run the construction of G\"{o}del's L instead of first-order logic, one obtains an inner model $C (\L^*)$ of $\mathsf{ZFC}$. 
Myhill and Scott showed that if $\L^*$ is full second-order logic, the inner model $C (\L^*)$ is equal to HOD, the class of all hereditarily ordinal definable sets (\cite{Myhill-Scott}). 
HOD is the largest inner model of $\mathsf{ZFC}$ such that its every element is ordinal definable. HOD is very complicated and very difficult to analyze in $\mathsf{ZFC}$.  
It is well-known that one can change HOD drastically in $\mathsf{ZFC}$ using set forcing, and that one cannot decide whether CH holds in HOD in $\mathsf{ZFC}$ even if one additionally assumes large cardinal axioms. 

Kennedy, Magidor, and V\"{a}\"{a}n\"{a}nen~\cite{MR4290501,part2} initiated a research program to investigate on the inner models $C(\L^*)$ for various logics $\L^*$ extending first-order logic. 
They have shown that $C(\L^*)$ is equal to G\"{o}del's L when $\L^*$ is $\mathcal{L} (Q_{\alpha})$ or $\mathcal{L} (Q_{\alpha}^{\text{MM}})$ (\cite{MR4290501}). They also introduced two interesting inner models $C^*$ (when $\L^*$ is $\mathcal{L}(Q^{\text{cf}}_{\omega})$) and $C (\text{aa})$ (when $\L^*$ is the stationary logic), and investigated on the basic properties of $C^*$ and $C (\text{aa})$ (\cite{MR4290501,part2}). 

In this paper, we introduce the inner model $C^{2b}$, that is $C(\L^*)$ when $\L^*$ is Boolean-valued second-order logic. We show that assuming that Projective Determinacy ($\mathsf{PD}$) holds in any set generic extension, the inner model $C^{2b}$ is the smallest inner model of $\mathsf{ZFC}$ closed under $\M_n^{\#}$ operators for all $n < \omega$. Also $C^{2b}$ is defined absolutely among all the set generic extensions. Furthermore, $C^{2b}$ admits a fine structure, hence it satisfies various combinatorial principles such as $\mathsf{GCH}$, $\diamondsuit_{\kappa}$ for all regular uncountable cardinals $\kappa$, $\square_{\kappa}$ for all uncountable cardinals $\kappa$ (Theorem~\ref{mainthm}). %Therefore, $C^{2b}$ enjoys various nice properties as G\"{o}del's L does while HOD enjoys none of them. 
%(or there are proper class many Woodin cardinals), the inner model $C^{2b}$ is the smallest inner model $N$ of $\mathsf{ZFC}$ such that for all $\mathbb{P} \in N$ and $\mathbb{P}$-generic filters $G$ over $V$, the structure $(H_{\omega_1}, \in)^{N[G]}$ is an elementary substructure of $(H_{\omega_1}, \in)^{V[G]}$ (or it is the smallest inner model $N$ of $\mathsf{ZFC}$ closed under $\M_n^{\#}$ operators for all $n < \omega$). In particular, $\mathsf{PD}$ holds in any set generic extension of $C^{2b}$. Also $C^{2b}$ is defined absolutely among all the set generic extensions. Furthermore, $C^{2b}$ admits a fine structure, hence it satisfies various combinatorial principles such as $\mathsf{GCH}$, $\diamondsuit_{\kappa}$ for all regular uncountable cardinals $\kappa$, $\square_{\kappa}$ for all uncountable cardinals $\kappa$ (Theorem~\ref{mainthm}). Therefore, $C^{2b}$ enjoys various nice properties as G\"{o}del's L does while HOD enjoys none of them. 

Let us summarize the difference between Boolean-valued second-order logic and full second-order logic in Table~\ref{table} below. 
\begin{table}[h]
%\begin{center}
\centering
\begin{threeparttable} 
\begin{tabular}{p{3cm}p{5cm}p{5cm}} \hline
    & Boolean-valued  \hspace{2cm} second-order logic & Full second-order logic  \\ \hline
  Validity  & consistently definable over the third-order arithmetic & $\Pi_2$-complete in $\mathsf{ZFC}$ \\ \hline
   Compactness number & \multicolumn{1}{c}{$\omega_1$} & the least extendible cardinal \\ \hline
   Hanf number & \multicolumn{1}{c}{$\omega$} & between the least supercompact cardinal and the least extendible cardinal\\ \hline 
   L\"{o}wenheim number & consistently below the least Woodin cardinal & the supremum of $\Delta_2$-definable ordinals in $\mathsf{ZFC}$\\ \hline
   Inner models & \multicolumn{1}{c}{$C^{2b}$} & \multicolumn{1}{c}{HOD} \\ \hline
\end{tabular}
%\vspace{0.5cm}
\caption{Boolean-valued second-order logic vs full second-order logic}\label{table}
%\end{center}
\end{threeparttable} 
\end{table}

All in all, we conclude that we have clarified various aspects and contexts where Boolean-valued second-order logic is much simpler than full second-order logic. 

This paper is organized as follows: In Section~\ref{sec:basics}, we introduce basic definitions and theorems we will use in the rest of the paper. 
In Section~\ref{sec:compactness}, we discuss the results on the compactness number of Boolean-valued second-order logic. 
In Section~\ref{sec:L2b}, we discuss the results on $C^{2b}$, the inner model constructed from Boolean-valued second-order logic. 
%In Section~\ref{sec:errors}, we mention two errors in the paper \lq Boolean-Valued Second-Order Logic~\cite{BVSOL}'. 
In Section~\ref{sec:Q}, we list some open questions.

%To be filled in.

\section{Basic definitions and theorems}\label{sec:basics}

From now on, we work in $\ZFC$ unless clearly specified. 

We assume that the reader is familiar with the basic theory of forcing (given in e.g.,~\cite{Jech}), basics of descriptive set theory (given in e.g.,~\cite{new_Moschovakis}), and basics of large cardinals (given in e.g.,~\cite{MR1994835}). For the proof of Proposition~\ref{prop:generically-extendible} in Section~\ref{sec:compactness}, we use the basic theory of stationary tower forcings (given in~\cite{MR2069032}). 
For the proofs of  Theorem~\ref{thm:PDmouse} and Theorem~\ref{thm:projective-absoluteness} in this section, and of Theorem~\ref{mainthm}~\ref{mainthm:item4} in Section~\ref{sec:L2b}, we assume the basics of inner model theory (given in e.g.,~\cite{outline-IMT} and \cite{FSIT}). 
%%% ADD PREREQUISITES AND BASIC REFERENCES

\begin{Def}[Boolean-valued structures]
Let $\mathcal{L} = \{R_i \mid i \in I\} \cup \{ c_j \mid j \in J\}$ be a language where each $R_i$ is a relation symbol and each $c_j$ is a constant symbol. A {\it Boolean-valued $\mathcal{L}$-structure} consists of $M = ( A, \mathbb{B}, \{ R_i^M \}_{i \in I}, \{ c_j^M\}_{j\in J})$ where
\begin{enumerate}
\item $A$ is a non-empty set, 

\item $B$ is a complete Boolean algebra, 

\item for each $i\in I$, if $R_i$ is an $m$-ary relation symbol, then $R_i^M $ is a function from $A^m$ to $\mathbb{B}$, and

\item for each $j \in J$, $c_j^M$ is a function from $A$ to $\mathbb{B}$ such that for all $a , b \in A$ with $a \neq b$, $c_j^M (a) \wedge c_j^M (b) = 0$ in $\mathbb{B}$, and that the set $c_j^M [A] \setminus \{ 0\}$ is a maximal antichain in the partial order $\mathbb{B} \setminus \{ 0\}$. 
\end{enumerate}
\end{Def}

When $\mathbb{B} = \{ 0, 1\}$, a Boolean-valued $\mathcal{L}$-structure can be seen as a first-order $\mathcal{L}$-structure by identifying $R_i^M$ with the characteristic function of an $m$-ary predicate and by identifying $c_j^M$ with the unique element $a$ of $A$ such that $c_j^M (a) = 1$. 

We now interpret second-order sentences by Boolean-valued structures in the following way:
\begin{Def}\label{Boolean-valued interpretation}
Let $\mathcal{L}$ be as above and $M = ( A, \mathbb{B} , \{R_i^M\}_{i \in I }, \{ c_j^M \}_{j\in J})$ be a Boolean-valued $\mathcal{L}$-structure. To each second-order $\mathcal{L}$-formula, $\vec{a} \in A^{< \omega}$, and $\vec{f} \in (\mathbb{B}^A)^{<\omega}$, we assign $\llbracket\phi [\vec{a}, \vec{f}]\rrbracket^M$ by induction on the complexity of $\phi$ in the following way:
\begin{enumerate}
\item $\phi$ is $x = y$ where $x,y$ are first-order variables. Then $\llbracket x= y [a,b]\rrbracket^M$ is $1$ if $a =b $ and $0$ otherwise.

\item $\phi$ is $c_j = x$ where $x$ is a first-order variable. Then $\llbracket c_j =x [a] \rrbracket^M = c_j^M (a)$. The case when $\phi$ is $x = c_j$ is similarly defined. 

\item $\phi$ is $c_{j_1} = c_{j_2}$. Then $\displaystyle \llbracket c_{j_1} = c_{j_2} \rrbracket^M = \bigvee_{a\in A} \bigl(c_{j_1}^M (a) \wedge c_{j_2}^M (a)\bigr)$. 

\item $\phi$ is $R_i (x)$ where $x$ is a first-order variable. Then $\llbracket R_i (x) [a]\rrbracket^M = R_i^M (a)$. 

\item $\phi$ is $X(x)$ where $X$ is a second-order variable and $x$ is a first-order variable. Then $\llbracket X(x) [a,f]\rrbracket^M = f(a)$. 

\item $\phi$ is $\neg \psi$. Then $\llbracket\phi [\vec{a},\vec{f}]\rrbracket^M = 1 - \llbracket \psi [\vec{a}, \vec{f}]\rrbracket^M$. 

\item $\phi$ is $\psi_1 \wedge \psi_2$. Then $\llbracket\phi [\vec{a},\vec{f}]\rrbracket^M = \llbracket \psi_1 [\vec{a}, \vec{f}]\rrbracket^M \wedge \llbracket\psi_2 [ \vec{a},\vec{f}]\rrbracket^M$. 

\item $\phi$ is $(\exists x) \, \psi$ where $x$ is a first-order variable. Then \[ \llbracket\phi [\vec{a},\vec{f}]\rrbracket^M = \bigvee_{b \in A} \llbracket\psi[b, \vec{a}, \vec{f}]\rrbracket^M.\] 

\item $\phi$ is $ (\exists X) \, \psi$ where $X$ is a second-order variable. Then \[ \llbracket\phi [\vec{a}, \vec{f}]\rrbracket^M = \bigvee_{g \colon A \to \mathbb{B}} \llbracket\psi [\vec{a}, g, \vec{f}]\rrbracket^M.\] 
\end{enumerate}
\end{Def}

Note that interpreting second-order formulas via Boolean-valued structures is the same as doing it via full second-order structures in set generic extensions:
\begin{Def}\label{def:MG}
Let $\mathcal{L}$ be as above and $M = ( A, \mathbb{B} , \{R_i^M\}_{i \in I }, \{ c_j^M \}_{j \in J})$ be a Boolean-valued $\mathcal{L}$-structure. Let $G$ be a $\mathbb{B}$-generic filter over $V$. Then we define the full second-order $\mathcal{L}$-structure $M[G]$ in $V[G]$ as follows:
\begin{enumerate}[(1)]
\item The first-order part of $M[G]$ is $A$, 

\item the second-order part of $M[G]$ is $\bigcup_{n \in \omega} \wp(A^n)^{V[G]}$,

\item the interpretation of $R_i$ is as follows:
\begin{align*}
R_i^{M[G]} = \{ a \in A \mid R_i^M (a) \in G\},
\end{align*}

\item the interpretation of $c_j$ is the unique element $a$ of $A$ with $c_j (a) \in G$.\label{MGitem4}
\end{enumerate}
\end{Def}
Notice that the existence of $a \in A$ with $c_j (a) \in G$ in Definition~\ref{def:MG}~\ref{MGitem4} follows from the condition that $c_j^M [A] \setminus \{ 0 \}$ is a maximal antichain in the partial order $\mathbb{B} \setminus \{ 0\}$ and the genericity of $G$, and that the uniqueness of $a \in A$ with $c_j (a) \in G$ in Definition~\ref{def:MG}~\ref{MGitem4} follows from the condition that for all $a, b \in A$ with $a \neq b$, $c_j^M (a) \wedge c_j^M (b) = 0$ in $\mathbb{B}$.

\begin{Lem}\label{boolean semantics}
Let $\mathcal{L}$ be as above and $M = ( A, \mathbb{B} , \{R_i^M\}_{i \in I }, \{ c_j^M\}_{j \in J})$ be a Boolean-valued $\mathcal{L}$-structure. Let $G$ be a $\mathbb{B}$-generic filter over $V$. Then for each second-order $\mathcal{L}$-formula, $\vec{a} \in A^{< \omega}$, and $\vec{f} \in (\mathbb{B}^A)^{<\omega}$, 
\begin{align*}
\llbracket\phi [\vec{a}, \vec{f}]\rrbracket^M \in G \iff M[G] \vDash \phi [\vec{a}, \vec{f^G}],
\end{align*}
where $\vec{f^G} = \{f_1^G, \ldots , f_n^G\}$ and for each $i$ with $1\le i \le n$,
\begin{align*}
f_i^G = \{ a \in A \mid f_i (a) \in G\}.
\end{align*}
\end{Lem}

\begin{proof}
It is simiar to \cite[Lemma~2.4]{BVSOL}.
\end{proof}

The following is an easy consequence of Lemma~\ref{boolean semantics} together with basics of forcing:
\begin{Lem}\label{generic extension to Boolean structure}
Let $\mathcal{L}$ be as above and $T$ be a set of second-order $\mathcal{L}$-sentences. Suppose that there is a set generic extension $V[G]$ of $V$ such that in $V[G]$, there is a full second-order $\mathcal{L}$-structure $N$ such that for all $\phi$ in $T$, $N \vDash \phi$. Then in $V$, there is a Boolean-valued model $M$ of $T$, i.e., $M$ is a Boolean-valued $\mathcal{L}$-structure such that $\bigwedge_{\phi  \in T} \, \llbracket\phi [\vec{a}, \vec{f}]\rrbracket^M  \neq 0$. 
\end{Lem}

\begin{proof}
Let $\mathcal{L}$, $T$, $G$, and $N$ be as in the statement. We may assume $G$ is $\mathbb{B}$-generic over $V$ for some complete Boolean algebra $\mathbb{B}$ in $V$. We may also assume that the first-order universe of $N$ is an ordinal $\gamma$. Let $\dot{N}$ be a $\mathbb{B}$-name with $\dot{N}^G = N$. Then in $V$, let $M = ( A , \mathbb{B} , \{R_i\}_{i \in I} , \{c_j \}_{j\in J})$ be be the following Boolean-valued $\mathcal{L}$-structure:
\begin{itemize}
\item $A$ is $\gamma$,

\item for each $i \in I$ and $a\in A$, let $R_i^M (a)$ be the Boolean value of the statement \lq $\check{a} \in R_i^{\dot{N}}$', and

\item for each $j \in J$ and $a \in A$, let $c_j^M (a)$ be the Boolean value of the statement \lq $\check{a} = c_j^{\dot{N}}$'.  
\end{itemize}
Then it follows that $M[G] = N$. Hence by Lemma~\ref{boolean semantics}, for all $\phi$ in $T$, $\llbracket\phi [\vec{a}, \vec{f}]\rrbracket^M  \in G$. By the genericity of $G$, it follows that $\bigwedge_{\phi  \in T} \, \llbracket\phi [\vec{a}, \vec{f}]\rrbracket^M  \neq 0$, as desired.  
\end{proof}

The following lemma will be used in the proofs of Theorem~\ref{thm:compactness2b} and Theorem~\ref{thm:L2bcompactness-number}:
\begin{Lem}\label{lem:Vbeta}
Let $\mathcal{L}$ be a language containing a binary predicate $E$. Then there is a second-order $\mathcal{L}$-sentence $\Phi_0$ such that if $P = (A, E^P , \ldots )$ is a full second-order $\mathcal{L}$-structure satisfying $\Phi_0$, then $(A, E^P)$ is isomorphic to $(V_{\beta} , \in)$ for some ordinal $\beta$.  
\end{Lem}

\begin{proof}
See e.g., the proof of \cite[Theorem~2]{compactness2ndordermagidor}.
\end{proof}

The following theorems will be used in Section~\ref{sec:L2b}. For the basics of $\M_n^{\#} (x)$ for an $n < \omega$ and a transitive set $x$, one can consult \cite[Section~7]{outline-IMT}, \cite{projective-IM}, and \cite[Section~2]{scale-KR}. 
\begin{Lem}\label{lem:Mnsharp}
%\begin{enumerate}[(1)]
\item Let $n < \omega$ and suppose that for all transitive sets $x$, $\M_n^{\#} (x)$ exists and it is fully iterable. 
Let $N$ be an inner model of $\mathsf{ZFC}$ which is closed under $\M_n^{\#}$ operator, i.e., for all transitive sets $x \in N$, we have $\M_n^{\#} (x) \in N$. 
Then the statement \lq\lq $y = \M_{n+1}^{\#} (x)$'' is absolute between $N$ and $V$.
%\item 
%\end{enumerate}
\end{Lem}

\begin{proof}
This is basically because the iteration strategy of $\M_{n+1}^{\#} (x)$ is guided by $\mathcal{Q}$-structures whose complexity is at most the $\M_n^{\#}$ operator. For the details, see \cite{outline-IMT} and \cite{projective-IM}. 
\end{proof}

\begin{Thm}[Martin; Harrington; Woodin; Neeman]\label{thm:PDmouse}

${}$

\begin{enumerate}[(1)]
\item Projective Determinacy ($\mathsf{PD}$) holds if and only if for all $n < \omega$ and all reals $x$, $\M_n^{\#} (x)$ exists and it is $\omega_1$-iterable.\label{PDitem1} 

\item The following are equivalent:\label{PDitem2}
\begin{enumerate}[(a)]
\item Projective Determinacy ($\mathsf{PD}$) holds in any set generic extension of $V$.\label{PDitema}

%\item For all set generic extensions $V[G]$ and $V$ and all set generic extensions $V[G][H]$ of $V[G]$, the structure $(H_{\omega_1},\in)^{V[G]}$ is an elementary substructure of $(H_{\omega_1}, \in)^{V[G][H]}$. 

\item In any set generic extension $V[G]$ of $V$, the following holds: for all $n \in \omega$ and all transitive sets $x$, $\text{M}^{\#}_n (x)$ exists and it is fully iterable.\label{PDitemb}
\end{enumerate}
\end{enumerate}
\end{Thm}

\begin{proof}
For \ref{PDitem1}, see \cite[Theorem~2.14]{optimal2} and \cite[Theorem~2.1]{MSW20}. 

For \ref{PDitem2}, we assume \ref{PDitema} and verify \ref{PDitemb}. By induction on $n < \omega$, we show that (i) for all set generic extensions $V[G]$ of $V$, the following is true in $V[G]$: for all transitive sets $x$, $\M_n^{\#} (x)$ exists and it is fully iterable, and (ii) for all set generic extensions $V[G]$ of $V$ and all set generic extensions $V[H]$ of $V[G]$, the statement \lq\lq $y = \M_n^{\#}(x)$'' is absolute between $V[G]$ and $V[G][H]$. 

Let $n = 0$. We verify (i) and (ii) above for $n =0$. Notice that (ii) for $n=0$ is the absoluteness of the statement \lq\lq $y = x^{\#}$'' between $V[G]$ and $V[G][H]$, which easily follows. %from the fact that a set forcing does not add a sharp of a given set. So we focus on verifying (i) for $n=0$. 
Let $G$ be as in the statement of (i). 
%Let $V[G]$ be any set generic extension of $V$ and $V[G][H]$ be any set generic extension of $V[G]$. 
We may assume $V[G] = V$ (the arguments for an arbitrary $V[G]$ is exactly the same as when $V[G] = V$).   We need to verify that for all transitive sets $x$, $x^{\#}$ exists and it is fully iterable. %(the absoluteness of the statement \lq\lq $y = x^{\#}$'' between $V$ and $V[H]$ easily follows). 
Let $x$ be a transitive set. Let $I$ be any $\text{Col} (\omega , x)$-generic filter over $V$. Then by assumption, $\mathsf{PD}$ holds in $V[I]$. So by \ref{PDitem1}, letting $a$ be a real coding $x$ in $V[I]$, we have $\M_0^{\#} (a)$ (which is $a^{\#}$) in $V[I]$. 
But $x^{\#}$ can be easily constructed from $a^{\#}$ in $V[I]$. So $x^{\#}$ exists in $V[I]$. Since a set forcing does not produce a sharp of a given set, $x^{\#}$ exists in $V$, and it is certainly fully iterable.  %and the statement \lq\lq $y = x^{\#}$'' is absolute between $V$ and $V[H]$. Hence $x^{\#}$ exists in $V$, and it is certainly fully iterable. 

Let $n < \omega$.   Assume (i) and (ii) above for $n$. We will verify (i) and (ii) above for $n+1$. First note that (ii) for $n+1$ follows from Lemma~\ref{lem:Mnsharp}, and (i) and (ii) for $n$. So we focus on verifying (i) for $n+1$. Let $V[G]$ be any set generic extension of $V$. We may assume $V[G] = V$ (the arguments for an arbitrary $V[G]$ is exactly the same as when $V[G] = V$). We will show that for all transitive sets $x$, $\M_{n+1}^{\#} (x)$ exists and it is fully iterable. Let $x$ be a transitive set. Let $I$ be any $\text{Col} (\omega , x)$-generic filter over $V$. Then by assumption, $\mathsf{PD}$ holds in $V[I]$. So by \ref{PDitem1}, letting $a$ be a real coding $x$ in $V[I]$, we have $\M_{n+1}^{\#} (a)$ in $V[I]$ and it is $\omega_1$-iterable in $V[I]$. But $\M_{n+1}^{\#} (x)$ can be easily constructed from $\M_{n+1}^{\#} (a)$ in $V[I]$. Furthermore, since the first-order theory of $\M_{n+1}^{\#} (x)$ in $V[I]$ (in the language of premice over $x$) does not depend on $I$ (but on the model $V[I]$), by the weak homogeneity of $\text{Col} (\omega , x)$, letting $M = \bigl(\M_{n+1}^{\#} (x)\bigr)^{V[I]}$, the first-order theory of $M$ (in the language of premice over $x$) belongs to $V$, and hence $M \in V$ because the first projectum of $M$ is equal to $x$.  
However, by (ii) for $n+1$, the statement \lq\lq $y = \M_{n+1}^{\#} (x)$'' is absolute between $V$ and $V[I]$. Hence $\M_{n+1}^{\#} (x)$ exists in $V$ and $M = \M_{n+1}^{\#} (x)$. 

Notice that the arguments in the last paragraph actually show that $M = \M_{n+1}^{\#} (x)$ is $\omega_1$-iterable in any set generic extension $V[I']$ of $V$. Furthermore, the iteration strategy of $M$ in $V[I']$ restricted to iteration trees in $V$ (which are countable in $V[I']$) must be in $V$ because it is guided by $\mathcal{Q}$-structures whose complexity is at most the $\M_n^{\#}$-operator, and $V$ is closed under the $\M_n^{\#}$ operator in $V[I']$ by (i) and (ii) for $n$. Hence $\M_{n+1}^{\#} (x)$ is fully iterable in $V$.  

This completes the verification of \ref{PDitemb} assuming \ref{PDitema}.

We now assume \ref{PDitemb} and verify \ref{PDitema}. Let $V[G]$ be any set generic extension of $V$. %via a partial order $\mathbb{P}$. 
By \ref{PDitem1}, it is enough to verify that in $V[G]$, for all $n < \omega$ and all reals $x$, $\M_n^{\#} (x)$ exists and it is $\omega_1$-iterable. But this dicretly follows from \ref{PDitemb}. 
%Fix any $n < \omega$ and any real $x$ in $V[G]$. Let $\tau$ be a $\mathbb{P}$-name with $\tau^G = x$. Let $a$ be the transitive closure of $\{ (\mathbb{P}, \tau) \}$. Then by \ref{PDitemb}, we have $\M_n^{\#} (a)$ in $V[G]$ and it is fully iterable. 
%Also, by \ref{PDitemb}, $V$ is closed under $\M_k^{\#}$ operators for all $k< \omega$ in $V[G]$. Hence by Lemma~\ref{lem:
%But in $V[G]$, one can construct $\M_n^{\#} (x)$ from $\M_n^{\#} (a)$ and $G$ in a simple manner, and $\M_n^{\#} (x)$ is fully iterable there. In particular, it is $\omega_1$-iterable, as desired. 

This completes the verification of \ref{PDitema} assuming \ref{PDitemb}.
\end{proof}

\begin{Rem}\label{rem:Mnsharp}
The above arguments for the direction from \ref{PDitema} to \ref{PDitemb} in \ref{PDitem2} in Theorem~\ref{thm:PDmouse} show the following: Suppose that $\mathsf{PD}$ holds in any set generic extension. Then for all set generic extensions $V[G]$ of $V$ %the $\M_n^{\#}$ operators are absolute between $V$ and $V[G]$ for all $n< \omega$, and so $V$ is closed under $\M_n^{\#}$ operators in $V[G]$ for all $n < \omega$.  
and all set generic extensions of $V[G][H]$ of $V[G]$, the $\M_n^{\#}$ operators are absolute between $V[G]$ and $V[G][H]$ for all $n< \omega$, and so $V[G]$ is closed under $\M_n^{\#}$ operators in $V[G][H]$ for all $n < \omega$.  
\end{Rem}

\begin{Thm}[Woodin]\label{thm:projective-absoluteness}
Suppose that $\mathsf{PD}$ holds in any set generic extension of $V$. %(or for all $n \in \omega$ and all transitive sets $x$, $\text{M}^{\#}_n (x)$ exists and is fully iterable). 
\begin{enumerate}[(1)]
\item Let $N$ be an inner model of $\mathsf{ZFC}$ such that for all $n < \omega$ and all reals $x \in N$, $\M_n^{\#} (x) \in N$. Then the structure $(H_{\omega_1}, \in)^N$ is an elementary substructure of $(H_{\omega_1}, \in)^V$.  \label{projective-absoluteness-item1}

\item Let $N$ be an inner model of $\mathsf{ZFC}$ which is closed under $\M_n^{\#}$ operators for all $n < \omega$, i.e., for all $n < \omega$ and all transitive sets $x \in N$, $\M_n^{\#} (x) \in N$. Then for all partial orders $\mathbb{P}$ in $N$ and all $\mathbb{P}$-generic filters $G$ over $V$, the structure $(H_{\omega_1},\in)^{N[G]}$ is an elementary substructure of $(H_{\omega_1},\in)^{V[G]}$. \label{projective-absoluteness-item2}

\item For all set generic extensions $V[G]$ and $V$ and all set generic extensions $V[G][H]$ of $V[G]$, the structure $(H_{\omega_1},\in)^{V[G]}$ is an elementary substructure of $(H_{\omega_1}, \in)^{V[G][H]}$. \label{projective-absoluteness-item3}
\end{enumerate}
\end{Thm}

\begin{proof}
For \ref{projective-absoluteness-item1}, see \cite[Section~7]{outline-IMT} and \cite{projective-IM}. 

For \ref{projective-absoluteness-item2}, let $N, \mathbb{P}, G$ be as in the statement. By \ref{projective-absoluteness-item1}, it is enough to verify that for all $n < \omega$ and all reals $x \in N[G]$, we have $\M_n^{\#} (x) \in N[G]$ in $V[G]$. Fix an $n < \omega$ and a real $x \in N[G]$. 
Let $\tau$ be a $\mathbb{P}$-name in $N$ with $\tau^G = x$ and let $a$ be the transitive closure of $\{ (\mathbb{P}, \tau)\}$. Then by assumption on $N$, $\M_n^{\#} (a) \in N$ in $V$. 
But in $V[G]$, one can construct $\M_n^{\#} (x)$ from $\M_n^{\#}(a)$ and $G$ in a simple manner so that $\M_n^{\#} (x) \in N[G]$, as desired. 

For \ref{projective-absoluteness-item3}, let $G$ and $H$ be as in the statement. By \ref{projective-absoluteness-item1}, it is enough to verify that $V[G]$ is closed under $\M_n^{\#}$ operators in $V[G][H]$ for all $n < \omega$. 
But this follows from Remark~\ref{rem:Mnsharp}.
\end{proof}

\begin{Thm}\label{thm:projective-uniformization}

${}$

\begin{enumerate}[(1)]
\item (Moschovakis) Suppose $\mathsf{PD}$ holds. Then every projective relation $R$ can be uniformized by a projective function $f$.\label{uniformization-item1}  %In particular, for each $n \in \omega$, the structure $(H_{\omega_1}, \in)$ has a $\Sigma_n$-Skolem function $f_n$ which is first-order definable over $(H_{\omega_1},\in)$.

\item Suppose that $\mathsf{PD}$ holds in any set generic extension of $V$. Then for every projective relation $R$, the uniformizing function $f$ for $R$ as in~\ref{uniformization-item1} above is invariant under set forcings, i.e., for all set generic extensions $V[G]$ of $V$ and for all $x \in H_{\omega_1}^V$, $f^V (x) = f^{V[G]} (x)$, and in $V[G]$, $f^{V[G]}$ uniformizes $R^{V[G]}$. \label{uniformization-item2} %for each $n \in \omega$, the $\Sigma_n$-Skolem function $f_n$ for $(H_{\omega_1}, \in)$ in~\ref{uniformization-item1} is invariant under set forcings, i.e., for all set generic extensions $V[G]$ of $V$ and for all $x \in H_{\omega_1}^V$, $f_n^V (x) = f_n^{V[G]} (x)$.
\end{enumerate}
\end{Thm}

\begin{proof}
For \ref{uniformization-item1}, see \cite[6C.5]{new_Moschovakis}. 

For \ref{uniformization-item2}, this follows from that the statements \lq\lq $f(x) = y$" and \lq\lq $f$ uniformizes $R$" are definable over $(H_{\omega_1}, \in)$, and that the structure $(H_{\omega_1}, \in)^V$ is an elementary substructure of $(H_{\omega_1},\in)^{V[G]}$ by Theorem~\ref{thm:projective-absoluteness}~\ref{projective-absoluteness-item3}. 
\end{proof}

%To be filled in.

\section{Compactness numbers}\label{sec:compactness}

In this section, we consider the numbers related to compactness of Boolean-valued second-order logic and of full second-order logic. 

For regular cardinals $\kappa , \lambda$, let $\mathrm{L}^2_{\kappa, \lambda}$ denote the second-order logic with $<$$\kappa$-long conjunctions and disjunctions and $<$$\lambda$-long sequences of (first-order and second-order) quantifiers when forming the set of formulas in this logic. The semantics of these formulas is defined in a standard way. 

Respectively, for regular cardinals $\kappa , \lambda$, let $\mathrm{L}^{2b}_{\kappa, \lambda}$ denote the Boolean-valued second-order logic with $<$$\kappa$-long conjunctions and disjunctions and $<$$\lambda$-long sequences of quantifiers when forming the set of formulas in this logic. The semantics of these formulas is the same as that of formulas in $\mathrm{L}^2_{\kappa, \lambda}$ except that of the second-order quantifiers. The semantics of the second-order quantifiers (or sequences of second-order quantifiers) is analogously given as in the definition of the semantics of Boolean-valued second-order logic. 

For a given logic $\mathrm{L}^*$ and an infinite cardinal $\kappa$, $\kappa$ is called {\it $\mathrm{L}^*$-compact} if for any set of sentences $T$ in $\mathrm{L}^*$, if any subset of $T$ with cardinality less than $\kappa$ has a model, then $T$ also has a model. The {\it compactness number of $\mathrm{L}^*$} is the least uncountable regular cardinal $\kappa$ such that $\kappa$ is $\mathrm{L}^*$-compact.

The following is a result due to Magidor on the compactness of $\mathrm{L}^{2}_{\kappa , \kappa}$ and full second-order logic $\mathrm{L}^2$ related to large cardinals:
\begin{Def}
Let $\kappa$ be an infinite cardinal.
\begin{enumerate}
\item Let $\alpha$ be an ordinal with $\kappa < \alpha$. We say $\kappa$ is {\it $\alpha$-extendible} if there are some ordinal $\beta$ and an elementary embedding $j \colon (V_{\alpha} , \in) \to (V_{\beta} , \in )$ with critical point $\kappa$ such that $j(\kappa) > \alpha$. 

\item We say $\kappa$ is {\it extendible} if for all ordinals $\alpha$ with $\kappa < \alpha$, $\kappa$ is $\alpha$-extendible.
\end{enumerate}
\end{Def}

\begin{Thm}[Magidor]\label{thm:Magidor-compactness}

${}$

\begin{enumerate}[(1)]
\item Let $\kappa$ be a regular uncountable cardinal. Then the following are equivalent:\label{Magidor1}
\begin{enumerate}
\item $\kappa$ is extendible,

\item $\kappa$ is $\mathrm{L}^2_{\kappa , \kappa}$-compact, and

\item $\kappa$ is $\mathrm{L}^2_{\kappa , \omega}$-compact.
\end{enumerate}
\item The compactness number of full second-order logic $\mathrm{L}^2$ is the least extendible cardinal if it exists.\label{Magidor2}
\end{enumerate}
\end{Thm}

\begin{proof}
See \cite[Theorem~4]{compactness2ndordermagidor}.
\end{proof}

We introduce the notion of {\it generically extendible cardinals}, and following the idea of Magidor in Theorem~\ref{thm:Magidor-compactness}, we obtain the corresponding result on the numbers related to compactness of $\mathrm{L}^{2b}_{\kappa , \kappa}$. %and Boolean-valued second-order logic $\mathrm{L}^{2b}$:
%following the idea of Magidor in Theorem~\ref{thm:Magidor-compactness}, we show that the compactness number of Boolean-valued second-order logic is the least generically extendible cardinal:
\begin{Def}
Let $\kappa$ be an infinite cardinal.
\begin{enumerate}
\item Let $\alpha$ be an ordinal with $\kappa < \alpha$. We say $\kappa$ is {\it generically $\alpha$-extendible} if there is a set generic extension $V[G]$ of $V$ such that in $V[G]$, there are an ordinal $\beta$ and an elementary embedding $j \colon (V_{\alpha} , \in)^V \to (V_{\beta} , \in)^{V[G]}$ with critical point $\kappa$ such that $j(\kappa) > \alpha$. 

\item We say $\kappa$ is {\it generically extendible} if for all ordinals $\alpha$ with $\kappa < \alpha$, $\kappa$ is generically $\alpha$-extendible.
\end{enumerate}
\end{Def}

\begin{Thm}\label{thm:compactness2b}
${}$

%\begin{enumerate}[(1)]
%\item 
Let $\kappa$ be a regular uncountable cardinal. Then the following are equivalent:%\label{compactness:item1}
\begin{enumerate}[(1)]
\item $\kappa$ is generically extendible,\label{compactness:itema}

\item $\kappa$ is $\mathrm{L}^{2b}_{\kappa , \kappa}$-compact, and\label{compactness:itemb}

\item $\kappa$ is $\mathrm{L}^{2b}_{\kappa , \omega}$-compact.\label{compactness:itemc}
\end{enumerate}
%\item The compatness number of Boolean-valued second-order logic $\mathrm{L}^{2b}$ exists if there is a generically extendible cardinal, and it is at most the least generically extendible cardinal.
%\end{enumerate}
\end{Thm}

\begin{proof}
The idea of the proof is the same as Magidor's for Theorem~\ref{thm:Magidor-compactness}~\ref{Magidor1}. 

%\begin{description}
%\item[\ref{compactness:item1}]

We show the directions from \ref{compactness:itema} to \ref{compactness:itemb}, from \ref{compactness:itemb} to \ref{compactness:itemc}, and from \ref{compactness:itemc} to \ref{compactness:itema}. 
\begin{description}
\item[\ref{compactness:itema}$\Rightarrow$\ref{compactness:itemb}]
   Let $\kappa$ be generically extendible. We show that $\kappa$ is $\mathrm{L}^{2b}_{\kappa , \kappa}$-compact.

   Let $\mathcal{L}$ be a %relational 
language and $T$ be a set of $\mathcal{L}$-sentences in $\mathrm{L}^{2b}_{\kappa , \kappa}$. Suppose that every subset $T_0$ of $T$ of size less than $\kappa$ has a Boolean-valued model. 
We argue that $T$ has a Boolean-valued model as well.

Let $\alpha$ be an ordinal with $\kappa < \alpha$ and $T \in V_{\alpha}$ such that $V_{\alpha}$ thinks every subset $T_0$ of $T$ of size less than $\kappa$ has a Boolean-valued model. Since $\kappa$ is generically extendible, in some generic extension $V[G]$ of $V$, there are an ordinal $\beta$ and an elementary embedding $j \colon (V_{\alpha}, \in)^V \to (V_{\beta}, \in)^{V[G]}$ with critical point $\kappa$ such that $j(\kappa) > \alpha$. 

By elementarity of $j$, $V_{\beta}^{V[G]}$ thinks every subset $T_0$ of $j(T)$ of size less than $j(\kappa)$ has a Boolean-valued model. Set $T_0 = j^{\lq\lq} T$. Then since $T \in V_{\alpha}$ and $j (\kappa) > \alpha$, the set $T_0 = j^{\lq\lq} T$ is of size less than $j(\kappa)$ in $V_{\beta}^{V[G]}$. So in $V_{\beta}^{V[G]}$, there is a Boolean-valued model $N$ of $j^{\lq\lq} T$, i.e., $N$ is a Boolean-valued $j( \mathcal{L} )$-structure of the form $N = (A , \mathbb{B}, \{ R_i^N \} )$ such that $\bigwedge_{\phi \in j^{\lq\lq} T} \llbracket \phi \rrbracket^N \neq 0$. 

Since $\bigwedge_{\phi \in j^{\lq\lq} T} \llbracket \phi \rrbracket^N \neq 0$ in $V[G]$, by Lemma~\ref{boolean semantics}, there is a $\mathbb{B}$-generic filter $H$ over $V[G]$ such that for all $\phi$ in $j^{\lq\lq} T$, $N[H] \vDash \phi$. Since $j$ is elementary, the critical point of $j$ is $\kappa$, and $j \in V[G]$, in $V[G][H]$, there is a full second-order $\mathcal{L}$-structure $P$ %of the form $\bigl(\delta , \bigcup_{n \in \omega} \wp (\gamma^n)^{V[G][H]}, \{ R_i^P \}\bigr)$ 
such that for all $\phi$ in $T$, $P \vDash \phi$. 
By Lemma~\ref{generic extension to Boolean structure}, in $V$, there is a Boolean-valued model of $T$, as desired. 

\item[\ref{compactness:itemb}$\Rightarrow$\ref{compactness:itemc}] This is clear.

\item[\ref{compactness:itemc}$\Rightarrow$\ref{compactness:itema}]
Let $\kappa$ be a regular uncountable cardinal which is $\mathrm{L}^{2b}_{\kappa , \omega}$-compact. We will show that $\kappa$ is generically extendigble.

Let $\alpha$ be an ordinal with $\kappa < \alpha$. We will argue that $\kappa$ is $\alpha$-generically extendible.

Let $\mathcal{L}$ be the language consisting of one binary predicate $E$ and the set of constanct symbols $\{ c_x \mid x \in V_{\alpha}\}$. 
Let $M_0 = \bigl(V_{\alpha} , \in , \{x \mid x \in V_{\alpha}\} \bigr)$ be the first-order $\mathcal{L}$-structure intepreting $E$ as $\in$ in $V_{\alpha}$ and $c_x$ as $x$ for all $x \in V_{\alpha}$. 
Let $EDiag$ be the elementary diagram of the first-order $\mathcal{L}$-structure $M_0$, i.e., the set of all first-order $\mathcal{L}$-sentences true in the structure $M_0$. %$\bigl(V_{\alpha_0} , \in , \kappa , \{x \mid x \in V_{\alpha_0}\} \bigr)$. 
Let $\mathcal{L}'$ be the language obtained by adding constant symbols $\{ d_{\xi} \mid \xi \le \alpha \}$) to $\mathcal{L}$, and $\Phi_0$ be the second-order $\mathcal{L}'$-sentence as in Lemma~\ref{lem:Vbeta}, i.e., if $P = (A, E^P , \ldots )$ is a full second-order $\mathcal{L}'$-structure satisfying $\Phi_0$, then $(A, E^P)$ is isomorphic to $(V_{\beta} , \in)$ for some ordinal $\beta$.

Now let $T = EDiag \cup \{ \Phi_0 \} \cup \{ d_{\xi} E d_{\eta} \mid \xi < \eta \le \alpha \} \cup \{ d_{\alpha} E c_{\kappa} \} \cup \{ (\forall x ) \ (x E c_{\xi} \to \bigvee_{\eta < \xi} x = c_{\eta}) \mid \xi < \kappa \}$. Then for all subsets $T_0$ of $T$ of size less than $\kappa$, $T_0$ has a full second-order model by expanding the first-order $\mathcal{L}$-structure $M_0 = \bigl(V_{\alpha} , \in , \{x \mid x \in V_{\alpha}\} \bigr)$ in a suitable way, namely, interpreting relavant $d_{\xi}$s to ordinals less than $\kappa$. In particular, $T_0$ has a Boolean-valued model (by letting the complete Boolean algebra be the trivial one $\{ 0, 1 \}$). 
Since $\kappa$ is $\mathrm{L}^{2b}_{\kappa , \omega}$-compact, there is a Boolean-valued model $M$ of $T$, i.e., $M = (A, \mathbb{B} , \ldots )$ is a Boolean-valued $\mathcal{L}'$-structure with $\bigwedge_{\phi \in T} \llbracket \phi \rrbracket^M \neq 0$. 

By Lemma~\ref{boolean semantics}, there is a $\mathbb{B}$-generic filter $G$ over $V$ such that for all $\phi$ in $T$, $M[G] \vDash \phi$. Since $\Phi_0 \in T$, $M[G] \vDash \Phi_0$, so by the property of $\Phi_0$, $(A, E^{M[G]})$ is isomorphic to $(V_{\beta}, \in)^{V[G]}$ for some ordinal $\beta$. 
We may assume that $(A, E^{M[G]}) = (V_{\beta} , \in)^{V[G]}$. 

Now let $j \colon V_{\alpha}^V \to V_{\beta}^{V[G]}$ be such that for all $x \in V_{\alpha}$, $j(x) = c_x^{M[G]}$. Then since $EDiag \subseteq T$, in $V[G]$, $j \colon (V_{\alpha} , \in)^V \to (V_{\beta}, \in)^{V[G]}$ is an elementary embedding. 

We claim that $\kappa$ is the critical point of $j$ and that $j(\kappa) > \alpha$, which would show that $\kappa$ is generically $\alpha$-extendible. %, and hence generically extendible. Let $\gamma = c^{M[G]}$. Then since $c_{\xi} E c$ is in $T$ for all $\xi < \alpha_0$, $\gamma$ is an ordinal at least $\alpha_0$. 
%On the other hand, since $K(c)$ is in $T$ and $K^{M_0} = \kappa$, 
Let $\gamma < \kappa$. Then since for each $\xi \le \gamma$, the sentence $(\forall x ) \ (x E c_{\xi} \to \bigvee_{\eta < \xi} x = c_{\eta})$ is in $T$, we have $j (\gamma) = c_{\gamma}^{M[G]} = \gamma$. 
On the other hand, since $\{ d_{\xi} E d_{\eta} \mid \xi < \eta \le \alpha \} \cup \{ d_{\alpha} E c_{\kappa} \} \subseteq T$, there must be at least $(\alpha + 1)$-many ordinals below $c_{\kappa}^{M[G]} = j(\kappa)$. Hence $j (\kappa) > \alpha$ and the critical point of $j$ is equal to $\kappa$, as desired. 
\end{description}
This completes the proof of Theorem~\ref{thm:compactness2b}. 
\end{proof}
%\end{description}

To investigate on the compactness number of Boolean-valued second-order logic $\mathrm{L}^{2b}$, we introduce a slightly weaker notion than generically extendible cardinals:
\begin{Def}
Let $\kappa$ be an infinite cardinal.
\begin{enumerate}
\item Let $\alpha$ be an ordinal with $\kappa < \alpha$. We say $\kappa$ is {\it weakly generically $\alpha$-extendible} if there is a set generic extension $V[G]$ of $V$ such that in $V[G]$, there are an ordinal $\beta$ and an elementary embedding $j \colon (V_{\alpha} , \in)^V \to (V_{\beta} , \in)^{V[G]}$ such that $j(\kappa) > \alpha$. (Notice that $\kappa$ may not be the critical point of $j$.)

\item We say $\kappa$ is {\it weakly generically extendible} if for all ordinals $\alpha$ with $\kappa < \alpha$, $\kappa$ is weakly generically $\alpha$-extendible.
\end{enumerate}
\end{Def}

Notice that if $\kappa$ is generically extendible, then so is weakly generically extendible. 
Also notice that if $\kappa$ and $\lambda$ are infinite cardinals with $\kappa < \lambda$ and $\lambda$ is weakly generically extendible, then $\kappa$ is weakly generically extendible as well. 
Hence the notion of weakly generically extendible cardinals is sensible only when we focus on the least weakly generically extendible cardinal (if it exists). 

\begin{Thm}\label{thm:L2bcompactness-number}
Let $\kappa$ be an infinite cardinal. Then $\kappa$ is weakly generically extendible if and only if $\kappa$ is $\mathrm{L}^{2b}$-compact. In particular, the compactness number of Boolean-valued second-order logic $\mathrm{L}^{2b}$ is the least weakly generically extendible cardinal (if it exists). 
\end{Thm}

\begin{proof}
%\begin{comment}
%It is enough to verify the following two statements:
%\begin{enumerate}[(1)]
%\item If $\kappa$ is a weakly generically extendible cardinal, then $\kappa$ is $\mathrm{L}^{2b}$-compact.\label{L2bcompactness1}
%
%\item If $\kappa$ is $\mathrm{L}^{2b}$-compact, then $\kappa$ is weakly generically extendible.\label{L2bcompactness2}
%\end{enumerate}

%We first verify \ref{L2bcompactness1}. 
The idea of the proof is the same as Theorem~\ref{thm:compactness2b}. 
The difference is that one does not need to consider infinite conjunctions and disjunctions of formulas while the price is that one cannot demand that the generic elementary embedding given by the compactness of $\mathrm{L}^{2b}$ has critical point $\kappa$. 

Let $\kappa$ be weakly generically extendible. We show that $\kappa$ is $\mathrm{L}^{2b}$-compact.

Let $\mathcal{L}$ be a %relational 
language and $T$ be a set of $\mathcal{L}$-sentences in second-order logic. Suppose that every subset $T_0$ of $T$ of size less than $\kappa$ has a Boolean-valued model. %, i.e., there is a Boolean-valued $\mathcal{L}$-structure $M$ % = ( A, \mathbb{B} , \{R_i^M\}_{1 \le i \le n })$ 
%such that $\bigwedge_{\phi \in T_0} \llbracket \phi \rrbracket^M \neq 0$. 
We argue that $T$ has a Boolean-valued model as well.%, in other words, there is some $\mathbb{B}$-generic filter $G$ over $V$ such that for all $\phi$ in $T_0$, $M[G] \vDash \phi$, by Lemma~\ref{boolean semantics}. 

Let $\alpha$ be an ordinal with $\kappa < \alpha$, $T \in V_{\alpha}$, and $|T| < \alpha$ such that $V_{\alpha}$ thinks every subset $T_0$ of $T$ of size less than $\kappa$ has a Boolean-valued model. Since $\kappa$ is weakly generically extendible, in some generic extension $V[G]$ of $V$, there are an ordinal $\beta$ and an elementary embedding $j \colon (V_{\alpha}, \in)^V \to (V_{\beta}, \in)^{V[G]}$ such that $j(\kappa) > \alpha$. 

By elementarity of $j$, $V_{\beta}^{V[G]}$ thinks every subset $T_0$ of $j(T)$ of size less than $j(\kappa)$ has a Boolean-valued model. Set $T_0 = j^{\lq\lq} T$. Then since $T \in V_{\alpha}$, $|T|< \alpha$, and $j (\kappa) > \alpha$, the set $T_0 = j^{\lq\lq} T$ is of size less than $j(\kappa)$ in $V_{\beta}^{V[G]}$. So in $V_{\beta}^{V[G]}$, there is a Boolean-valued model $N$ of $j^{\lq\lq} T$, i.e., $N$ is a Boolean-valued $j( \mathcal{L} )$-structure of the form $N = (A , \mathbb{B}, \{ R_i^N \} )$ such that $\bigwedge_{\phi \in j^{\lq\lq} T} \llbracket \phi \rrbracket^N \neq 0$. %We may assume that the first-order universe of $N$ is an ordinal, so $N$ is of the form $(\delta , \mathbb{B}, \{ R_i^N \})$. 

Since $\bigwedge_{\phi \in j^{\lq\lq} T} \llbracket \phi \rrbracket^N \neq 0$ in $V[G]$, by Lemma~\ref{boolean semantics}, there is a $\mathbb{B}$-generic filter $H$ over $V[G]$ such that for all $\phi$ in $j^{\lq\lq} T$, $N[H] \vDash \phi$. Since $j$ is elementary and $j \in V[G]$, in $V[G][H]$, there is a full second-order $\mathcal{L}$-structure $P$ %of the form $\bigl(\delta , \bigcup_{n \in \omega} \wp (\gamma^n)^{V[G][H]}, \{ R_i^P \}\bigr)$ 
such that for all $\phi$ in $T$, $P \vDash \phi$. 
By Lemma~\ref{generic extension to Boolean structure}, in $V$, there is a Boolean-valued model of $T$, as desired. 
%To be filled in.

Let $\kappa$ be an infinite cardinal which is $\mathrm{L}^{2b}$-compact. We show that $\kappa$ is weakly generically extendible. 

Let $\alpha$ be an ordinal with $\kappa < \alpha$. We will argue that $\kappa$ is weakly $\alpha$-generically extendible.

%Let $\alpha_0$ be an ordinal with $\kappa < \alpha_0$ such that for all $\mu < \alpha$, if $\mu$ is not generically extendible, then there is some $\alpha < \alpha_0$ such that $\mu$ is not generically $\alpha$-extendible. We will argue that there is a $\mu \le \kappa$ such that $\mu$ is generically $\alpha_0$-extendible. By the property of $\alpha_0$, this would give us that $\mu$ is generically extendible. 

As is in the proof of \ref{compactness:itemc}$\Rightarrow$\ref{compactness:itema} in Theorem~\ref{thm:compactness2b}, 
let $\mathcal{L}$ be the language consisting of one binary predicate $E$ and the set of constanct symbols $\{ c_x \mid x \in V_{\alpha}\}$. 
Let $M_0 = \bigl(V_{\alpha} , \in , \{x \mid x \in V_{\alpha}\} \bigr)$ be the first-order $\mathcal{L}$-structure intepreting $E$ as $\in$ in $V_{\alpha}$ and $c_x$ as $x$ for all $x \in V_{\alpha}$. 
Let $EDiag$ be the elementary diagram of the first-order $\mathcal{L}$-structure $M_0$, i.e., the set of all first-order $\mathcal{L}$-sentences true in the structure $M_0$. %$\bigl(V_{\alpha_0} , \in , \kappa , \{x \mid x \in V_{\alpha_0}\} \bigr)$. 
Let $\mathcal{L}'$ be the language obtained by adding constant symbols $\{ d_{\xi} \mid \xi \le \alpha \}$) to $\mathcal{L}$, and $\Phi_0$ be the second-order $\mathcal{L}'$-sentence as in Lemma~\ref{lem:Vbeta}, i.e., if $P = (A, E^P , \ldots )$ is a full second-order $\mathcal{L}'$-structure satisfying $\Phi_0$, then $(A, E^P)$ is isomorphic to $(V_{\beta} , \in)$ for some ordinal $\beta$.

%Let $\mathcal{L}$ be the language consisting of one binary predicate $E$, one unary predicate $K$, and the set of constanct symbols $\{ c_x \mid x \in V_{\alpha_0}\}$. 
%Let $\bigl(V_{\alpha_0} , \in , \kappa , \{x \mid x \in V_{\alpha_0}\} \bigr)$ be the first-order $\mathcal{L}$-structure intepreting $E$ as $\in$ in $V_{\alpha_0}$, $K$ as $\kappa$ (as a set of ordinals), $c_x$ as $x$ for all $x \in V_{\alpha_0}$. 
%Let $EDiag$ be the elementary diagram of the first-order $\mathcal{L}$-structure $M_0 = \bigl(V_{\alpha_0} , \in , \kappa , \{x \mid x \in V_{\alpha_0}\} \bigr)$, i.e., the set of all first-order $\mathcal{L}$-sentences true in the structure $\bigl(V_{\alpha_0} , \in , \kappa , \{x \mid x \in V_{\alpha_0}\} \bigr)$. 
%Let $\mathcal{L}'$ be the language obtained by adding a constant symbol $c$ to $\mathcal{L}$, and $\Phi_0$ be the second-order $\mathcal{L}'$-sentence as in Lemma~\ref{lem:Vbeta}, i.e., if $P = (A, E^P , \ldots )$ is a full second-order $\mathcal{L}'$-structure satisfying $\Phi_0$, then $(A, E^P)$ is isomorphic to $(V_{\beta} , \in)$ for some ordinal $\beta$.

Now let $T = EDiag \cup \{ \Phi_0 \} \cup \{ d_{\xi} E d_{\eta} \mid \xi < \eta \le \alpha \} \cup \{ d_{\alpha} E c_{\kappa} \}$. %\cup \{ (\forall x ) \ (x E c_{\xi} \to \bigvee_{\eta < \xi} x = c_{\eta}) \mid \xi < \kappa \}$. 
Then for all subsets $T_0$ of $T$ of size less than $\kappa$, $T_0$ has a full second-order model by expanding the first-order $\mathcal{L}$-structure $M_0 = \bigl(V_{\alpha} , \in , \{x \mid x \in V_{\alpha}\} \bigr)$ in a suitable way, namely, interpreting relavant $d_{\xi}$s to ordinals less than $\kappa$. In particular, $T_0$ has a Boolean-valued model (by letting the complete Boolean algebra be the trivial one $\{ 0, 1 \}$). 
Since $\kappa$ is $\mathrm{L}^{2b}$-compact, there is a Boolean-valued model $M$ of $T$, i.e., $M = (A, \mathbb{B} , \ldots )$ is a Boolean-valued $\mathcal{L}'$-structure with $\bigwedge_{\phi \in T} \llbracket \phi \rrbracket^M \neq 0$. 

%Now let $T = EDiag \cup \{ \Phi_0 \} \cup \{ c_{\xi} E c \mid \xi < \kappa \} \cup \{ K(c) \}$. Then for all subsets $T_0$ of $T$ of size less than $\kappa$, $T_0$ has a full second-order model by expanding the first-order $\mathcal{L}$-structure $\bigl(V_{\alpha_0} , \in , \kappa , \{x \mid x \in V_{\alpha_0}\} \bigr)$ in a suitable way. In particular, $T_0$ has a Boolean-valued model (by letting the complete Boolean algebra be the trivial one $\{ 0, 1 \}$). 
%Since $\kappa$ is $\mathrm{L}^{2b}$-compact, there is a Boolean-valued model $M$ of $T$, i.e., $M = (A, \mathbb{B} , \ldots )$ is a Boolean-valued $\mathcal{L}'$-structure with $\bigwedge_{\phi \in T} \llbracket \phi \rrbracket^M \neq 0$. 

By Lemma~\ref{boolean semantics}, there is a $\mathbb{B}$-generic filter $G$ over $V$ such that for all $\phi$ in $T$, $M[G] \vDash \phi$. Since $\Phi_0 \in T$ and $M[G] \vDash \Phi_0$,  by the property of $\Phi_0$, $(A, E^{M[G]})$ is isomorphic to $(V_{\beta}, \in)^{V[G]}$ for some ordinal $\beta$. 
We may assume that $(A, E^{M[G]}) = (V_{\beta} , \in)^{V[G]}$. 

Now let $j \colon V_{\alpha_0}^V \to V_{\beta}^{V[G]}$ be such that for all $x \in V_{\alpha_0}$, $j(x) = c_x^{M[G]}$. Then since $EDiag \subseteq T$, in $V[G]$, $j \colon (V_{\alpha} , \in)^V \to (V_{\beta}, \in)^{V[G]}$ is an elementary embedding. 

We claim that $j(\kappa) > \alpha$, which would show that $\kappa$ is weakly generically $\alpha$-extendible. %, and hence generically extendible. Let $\gamma = c^{M[G]}$. Then since $c_{\xi} E c$ is in $T$ for all $\xi < \alpha_0$, $\gamma$ is an ordinal at least $\alpha_0$. 
%On the other hand, since $K(c)$ is in $T$ and $K^{M_0} = \kappa$, 
Since $\{ d_{\xi} E d_{\eta} \mid \xi < \eta \le \alpha \} \cup \{ d_{\alpha} E c_{\kappa} \} \subseteq T$, there must be at least $(\alpha + 1)$-many ordinals below $c_{\kappa}^{M[G]} = j(\kappa)$. Hence $j (\kappa) > \alpha$, as desired. 
%Let $\delta$ be the critical point of $j$. We claim that $\delta \le \kappa$ and $j(\kappa) > \alpha_0$, which would show that $\delta$ is generically $\alpha_0$-extendible, and hence generically extendible. Let $\gamma = c^{M[G]}$. Then since $c_{\xi} E c$ is in $T$ for all $\xi < \alpha_0$, $\gamma$ is an ordinal at least $\alpha_0$. 
%On the other hand, since $K(c)$ is in $T$ and $K^{M_0} = \kappa$, 
%\end{comment}
\end{proof}

We do not know whether the existence of weakly generically extendible cardinals is equivalent to that of generically extendible cardinals. 
However, we will show that they are equiconsistent as in Corollary~\ref{cor:equiconsistency}.

While extendible cardinals are very strong large cardinals, generically extendible cardinals can be $\omega_1$ under the presence of large cardinals:
\begin{Prop}\label{prop:generically-extendible}
Suppose there are proper class many Woodin cardinals. Then {\it every} regular uncountable cardinal is generically extendible.
\end{Prop}

\begin{proof}
We use the basic theory of stationary tower forcings given in~\cite{MR2069032}.

Let $\kappa$ be a regular uncountable cardinal. We show that $\kappa$ is generically extendible.
Let $\alpha$ be an ordinal with $\kappa < \alpha$. We will argue that $\kappa$ is $\alpha$-generically extendible.

Let $\delta$ be a Woodin cardinal with $\alpha <\delta$ and $\mathbb{P}_{<\delta}$ be the stationary tower forcing whose conditions $S$ in $V_{\delta}$ are stationary sets in $\bigcup S$. %such that $S$ consists of $Z$ such that $|Z | = \lambda \subseteq Z$. 
Let $S_0 = \{ Z \subseteq \alpha \mid Z \cap \kappa \in \kappa \text{ and the order type of $Z$ is less than $\kappa$} \}$. Then since $\kappa$ is regular uncountable, $S_0$ is stationary in $\alpha$ and so $S_0 \in \mathbb{P}_{<\delta}$.

Let $G$ be a $\mathbb{P}_{<\delta}$-generic filter over $V$ with $S_0 \in G$. Let $j \colon V \to M \subseteq V[G]$ be the generic elementary embedding induced by $G$. 
Then $j(\delta) = \delta$ and $M$ is closed under $<$$\delta$-sequences in $V[G]$, so $V_{\delta}^M = V_{\delta}^{V[G]}$. 

Since $S_0 \in G$, $j^{\lq\lq} \bigcup S_0 \in j (S_0)$. Since $\bigcup S_0 = \alpha$, $j^{\lq\lq}\alpha \in j(S_0)$, so $j^{\lq\lq} \alpha \cap j(\kappa) \in j(\kappa)$, and the order type of $j^{\lq\lq} \alpha$ is less than $j(\kappa)$. Since $j^{\lq\lq} \alpha \cap j(\kappa) = j^{\lq\lq} \kappa$, which is an element of $j(\kappa)$, $ j^{\lq\lq} \kappa$ is an ordinal, so  $j\upharpoonright \kappa = \text{id}$. Since the order type of $j^{\lq\lq} \alpha$ is $\alpha$, we have $j(\kappa) > \alpha$. So the critical point of $j$ is $\kappa$. 

Let $\beta = j(\alpha)$ and $k = j \upharpoonright V_{\alpha}$. Then $k$ is in $V[G]$, $k \colon (V_{\alpha}, \in)^V \to (V_{\beta}, \in)^M$ is elementary, the critical point of $k$ is $\kappa$, and $k(\kappa) > \alpha$. Since $\alpha < \delta$ and $j(\delta) = \delta$, we have $\beta = k (\alpha) = j(\alpha) < \delta$. Since $V_{\delta}^M = V_{\delta}^{V[G]}$, $V_{\beta}^M = V_{\beta}^{V[G]}$. Hence $k \colon  (V_{\alpha}, \in)^V \to (V_{\beta}, \in)^{V[G]}$ is an elementary embedding in $V[G]$ such that the critical point of $k$ is $\kappa$ and $k(\kappa) > \alpha$. Therefore, $\kappa$ is $\alpha$-generically extendible, as desired.
\end{proof}

\begin{Cor}\label{cor:compactness-BVSOL}
Suppose that there are proper class many Woodin cardinals. Then the compactness number of Boolean-valued second-order logic $\mathrm{L}^{2b}$ exists and it is equal to $\omega_1$.
\end{Cor}

\begin{proof}
By Theorem~\ref{thm:L2bcompactness-number}, the compactness number of Boolean-valued second-order logic $\mathrm{L}^{2b}$ is the least weakly generically extendible cardinal if it exists.
By Proposition~\ref{prop:generically-extendible}, $\omega_1$ is the least generically extendible cardinal. Hence the compactness number of Boolean-valued second-order logic $\mathrm{L}^{2b}$ exists and it is at most $\omega_1$. Since any generic elementary embedding of the form $j \colon (V_{\alpha},\in)^V \to (V_{\beta},\in)^{V[G]}$ fixes $\omega$, i.e., $j(\omega ) = \omega$, weakly generically extendible cardinals cannot be $\omega$. Hence the compactness number of Boolean-valued second-order logic $\mathrm{L}^{2b}$ is equal to $\omega_1$, as desired.
%Hence it is enough to verify that $\omega$ cannot be $\mathrm{L}^{2b}$-compact.
\end{proof}

We now discuss the consistency strength of the existence of the compactness number of Boolean-valued second-order logic $\mathrm{L}^{2b}$. 

First, Usuba~\cite{Usuba} proved that the consistency strength of the existence of generically extendible cardinals is the same as that of {\it virtually extendible cardinals}:
\begin{Def}[Bagaria-Gitman-Schindler~\cite{MR3598793}, Gitman-Schindler~\cite{MR3860539}]
Let $\kappa$ be an infinite cardinal. We say $\kappa$ is {\it virtually extendible} if for all ordinals $\alpha$ with $\kappa < \alpha$, there is a set generic extension $V[G]$ of $V$ where the following holds: There are an ordinal $\beta$ and an elementary embedding $j \colon (V_{\alpha} , \in )^V \to (V_{\beta}, \in)^V$ with critical point $\kappa$ such that $j(\kappa) > \alpha$.  
\end{Def}

\begin{Thm}[Usuba]\label{thm:Usuba}
The following are equiconsistent:
\begin{enumerate}
\item $\mathsf{ZFC}$+\lq\lq There is a virtually extendible cardinal''.

\item $\mathsf{ZFC}$+\lq\lq There is a generically extendible cardinal''.

\item $\mathsf{ZFC}$+\lq\lq $\omega_1$ is generically extendible''.

\item $\mathsf{ZFC}$+\lq\lq $\omega_2$ is generically extendible''.
\end{enumerate}
\end{Thm}

\begin{proof}
See \cite{Usuba}. %\cite[Theorem~1.4]{Usuba}.
\end{proof}

In \cite{MR3598793} and \cite{MR3860539}, they proved that if $\kappa$ is virtually extendible, then so is in $\mathrm{L}$. So the existence of a virtually extendible cardinal does not imply $0^{\#}$. 
Moreover, they proved that every Silver indiscernible is virtually extendible.

\begin{Cor}\label{L2bcompactness-zerosharp}
The existence of the compactness number of Boolean-valued second-order logic does not imply $0^{\#}$. 
\end{Cor}

\begin{proof}
This follows from Theorem~\ref{thm:L2bcompactness-number}, Theorem~\ref{thm:Usuba}, and the paragraph right before this corollary.
\end{proof}

We do not know if the existence of weakly generically extendible cardinals (equivalently, that of the compactness number of Boolean-valued second-order logic $\mathrm{L}^{2b}$) implies that of generically extendible cardinals. 
However, the following result shows that they are equiconsistent. 
The author would like to thank Toshimichi Usuba for communicating with him on the following result:
\begin{Prop}[Usuba]\label{prop:Usuba}
Suppose that there is a weakly generically extendible cardinal while there is no generically extendible cardinal. Then $0^{\#}$ exists.
\end{Prop}

For the proof of Proposition~\ref{prop:Usuba}, we use the following theorem of Kunen:
\begin{Thm}[Kunen]\label{thm:Kunen}
The following are equivalent:
\begin{enumerate}
\item $0^{\#}$ exists.

\item There are ordinals $\alpha , \beta$, and an elementary embedding $j \colon (\mathrm{L}_{\alpha}, \in) \to (\mathrm{L}_{\beta}, \in)$ such that the critical point of $j$ is less than $|\alpha|$, the cardinality of $\alpha$.
\end{enumerate}
\end{Thm}

\begin{proof}[Proof of Theorem~\ref{thm:Kunen}]
See e.g., \cite[21.1~Theorem]{MR1994835}.
\end{proof}

\begin{proof}[Proof of Proposition~\ref{prop:Usuba}]
By Theorem~\ref{thm:Kunen}, it is enough to show that there is a set generic extension $V[G]$ of $V$ where there are ordinals $\alpha , \beta$ and an elementary embedding $j \colon (\mathrm{L}_{\alpha} , \in) \to (\mathrm{L}_{\beta}, \in)$ such that the critical point of $j$ is less than the cardinality of $\alpha$ in $V[G]$ (then $0^{\#}$ would exist in $V[G]$ and hence in $V$). 

Assume that there is a weakly generically extendible cardinal $\kappa$ while there is no generically extendible cardinal. 
%Let $\kappa$ be the least weakly generically extendible cardinal. 
Since there is no generically extendible cardinal, there is an ordinal $\alpha$ with $\kappa < \alpha$ such that for all $\xi \le \kappa$, $\xi$ is {\it not} generically $\alpha$-extendible. 

%Let $\alpha$ be an ordinal with $\gamma < \alpha$.
Since $\kappa$ is weakly generically extendible, there is a set generic extension $V[G]$ of $V$ where there are an ordinal $\beta$ and an elementary embedding $j \colon (V_{\alpha}, \in)^V \to (V_{\beta} , \in)^{V[G]}$ such that $j(\kappa) > \alpha$. 
Since $\kappa$ is {\it not} generically $\alpha$-extendible, the critical point of $j$ cannot be equal to $\kappa$. Hence the critical point of $j$ must be smaller than $\kappa$.

Let $\mu$ be the critical point of $j$. Then $\mu < \kappa$. By the property of $\alpha$, $\mu$ is {\it not} generically $\alpha$-extendible, so $j(\mu) \le   \alpha$. 
Since $\mu$ is a cardinal in $V_{\alpha}^V$, by elementarity of $j$, $j(\mu)$ is a cardinal in $V_{\beta}^{V[G]}$, and hence in $V[G]$. But $j(\mu) \le \alpha$. So $\mu$ is less than the cardinality of $\alpha$ in $V[G]$.

Let $k = j \upharpoonright \mathrm{L}_{\alpha}$. Then in $V[G]$, $k \colon (\mathrm{L}_{\alpha}, \in) \to (\mathrm{L}_{\beta}, \in)$ is an elementary embedding such that the critical point of $k$ (which is $\mu$) is less than the cardinality of $\alpha$, as desired.
%To be filled in.
\end{proof}

\begin{Cor}\label{cor:equiconsistency}
The following are equiconsistent:
\begin{enumerate}
\item $\mathsf{ZFC}$+\lq\lq There is a generically extendible cardinal''.

\item $\mathsf{ZFC}$+\lq\lq There is a weakly generically extendible cardinal''.

\item $\mathsf{ZFC}$+\lq\lq The compactness number of Boolean-valued second-order logic $\mathrm{L}^{2b}$ exists''.
\end{enumerate}
\end{Cor}

\begin{proof}
This follows from Theorem~\ref{thm:L2bcompactness-number}, Theorem~\ref{thm:Usuba}, Corollary~\ref{L2bcompactness-zerosharp}, and Proposition~\ref{prop:Usuba}.
\end{proof}

In summary of this section, we have seen the following:
\begin{enumerate}
\item The compactness number of Boolean-valued second-order logic $\mathrm{L}^{2b}$ is equal to $\omega_1$ under the presence of proper class many Woodin cardinals while the compactness number of full second-order logic $\mathrm{L}^2$ is the least extendible cardinal.

\item The existence of the compactness number of Boolean-valued second-order logic $\mathrm{L}^{2b}$ is equiconsistent with that of generically extendible cardinals, and neither of them implies $0^{\#}$ while the existence of the compactness number of full second-order logic $\mathrm{L}^2$ is equiconsistent with that of extendible cardinals, which is much stronger than the existence of $0^{\#}$. 
\end{enumerate}

\section{The construction of G\"{o}del's Constructible Universe}\label{sec:L2b}

In this section, we discuss the inner model $C^{2b}$ obtained from the construction of G\"{o}del's Constructible Universe using Boolean-valued second-order logic $\mathrm{L}^{2b}$ instead of first-order logic. 

\begin{Def}[Kennedy, Magidor, and V\"{a}\"{a}n\"{a}nen~\cite{MR4290501,part2}]\label{def:L-hierarchy}
Let $\mathrm{L}^*$ be a logic. For a set $M$, let $\text{Def}_{\mathrm{L}^*} (M)$ be the set of all sets of the form $X = \{ a \in M \mid (M,\in) \vDash \phi [a, \vec{b} ] \}$, where $\phi (x, \vec{y})$ is an arbitrary formula of the logic $\mathrm{L}^*$ and $\vec{b} \in M^{<\omega}$. 
We define the hierarchy $(\mathrm{L}'_{\alpha} \mid \alpha \in \text{Ord})$ as follows:
\begin{align*}
\mathrm{L}'_0 =& \ \emptyset \\
\mathrm{L}'_{\alpha+1} = & \ \text{Def}_{\mathrm{L}^*} (\mathrm{L}'_{\alpha}) \\
\mathrm{L}'_{\gamma} = & \bigcup_{\alpha < \gamma} \mathrm{L}'_{\alpha} \text{ if $\gamma$ is a limit ordinal}.
\end{align*}
We write $C_o (\mathrm{L}^*)$ to denote the class $\bigcup_{\alpha \in \text{Ord}} \mathrm{L}'_{\alpha}$.\footnote{In \cite{part2}, Kennedy, Magidor, and V\"{a}\"{a}n\"{a}nen introduced a different inner model $C (\mathrm{L}^*)$ using Jensen's hierarchy and its truth predicate to make sure that the inner model satisfies the Axiom of Choice. For the logics $\L^*$ we consider in this paper, the inner models $C_o (\mathrm{L}^*)$ and $C (\mathrm{L}^*)$ coincide and both of them satisfy the Axiom of Choice.}
\end{Def}

If $\mathrm{L}^*$ is first-order logic, then the inner model $C_o (\mathrm{L}^*)$ is the G\"{o}del's Constructible Universe $\mathrm{L}$. 
Myhill and Scott~\cite{Myhill-Scott} proved that if $\mathrm{L}^*$ is full second-order logic $\mathrm{L}^2$, then the inner model $C_o (\mathrm{L}^*)$ is equal to HOD, the class of all sets which are hereditarily ordinal definable.

To discuss the inner model $C_o (\mathrm{L}^{2b})$ where $\mathrm{L}^{2b}$ is Boolean-valued second-order logic, we introduce the notion of definability in Boolean-valued second-order logic:
\begin{Def}
Let $\mathcal{L} = \{ E \}$ for a binary predicate $E$, and $A$ be a set.
\begin{enumerate}
\item Let $\vec{a} \in A^{<\omega}$ and $\phi$ be a second-order $\mathcal{L}$-formula. We say $\phi [\vec{a}]$ is {\it Boolean-valid with respect to $(A,\in)$} if for all Boolean-valued $\mathcal{L}$-structures $M$ with first-order structure $(A, \in)$ (i.e., the first-order universe of $M$ is $A$ and for all $y,z \in A$, $E^M (y,z) = 1$ if $y \in z$ and $E^M (y,z) = 0$ if $y\notin z$), we have $\Qp{\phi [\vec{a}]}^M = 1$. 

\item  Let $\vec{a} \in A^{<\omega}$ and $\phi$ be a second-order $\mathcal{L}$-formula. We say $(\phi , \vec{a})$ is {\it suitable to $(A,\in)$} if for all $x \in A$, either $\phi [x,\vec{a}]$ is Boolean-valid with respect to $(A,\in)$ or $\neg \phi [x, \vec{a}]$ is Boolean-valid with respect to $(A,\in)$.

\item Let $\text{Def}_{\mathrm{L}^{2b}} (A)$ be the set of sets of the form $X = \{ x \in A \mid \phi [x, \vec{a}]$ is Boolean-valid with respect to $(A,\in)\}$ for some $(\phi ,\vec{a})$ suitable to $(A,\in)$. 

\item For ordinals $\alpha$, we write $\mathrm{L}^{2b}_{\alpha}$ for $\mathrm{L}'_{\alpha}$ in Definition~\ref{def:L-hierarchy} when $\mathrm{L}^*$ is $\mathrm{L}^{2b}$. We also write $C^{2b}$ for the class $C_o (\mathrm{L}^*)$ in Definition~\ref{def:L-hierarchy} when $\mathrm{L}^*$ is $\mathrm{L}^{2b}$.%in this case.

\item For a set generic extension $V[G]$ of $V$, we write $N^{V[G]}_A$ for the full second-order $\mathcal{L}$-structure $(A, E^{\in} , \wp (A) , \in )^{V[G]}$ in $V[G]$, where $E^{\in} = \{ (y ,z ) \in A^2 \mid y \in z \}$.  
\end{enumerate}
\end{Def}

\begin{Rem}\label{rem:boolean-validity}
Notice that by Lemma~\ref{boolean semantics}, $\phi [\vec{a}]$ is Boolean-valid with respect to $(A,\in)$ if and only if for all set generic extensions $V[G]$ of $V$, we have $N^{V[G]}_A \vDash \phi [\vec{a}]$. %in $V[G]$, for all full second-order $\mathcal{L}$-structures $N$ with first-order universe $A$ and $E^N = \{ (y,z) \in A^2 \mid y \in z\}$, we have $N \vDash \phi [\vec{a}]$.
\end{Rem}   

Notice also that for all sets $A$, the family $\text{Def}_{\mathrm{L}^{2b}} (A)$ contains all the first-order definable subsets of $(A,\in)$. Hence the class $C^{2b}$ is an inner model of $\mathsf{ZF}$. 

We now state the main theorem of this section:
\begin{Thm}\label{mainthm}
Suppose that $\mathsf{PD}$ holds in any set generic extension of $V$. Then the following hold: 
\begin{enumerate}[(1)]
\item For all set generic extensions $V[G]$ of $V$ and all ordinals $\alpha$, $\bigl(\mathrm{L}^{2b}_{\alpha}\bigr)^V = \bigl(\mathrm{L}^{2b}_{\alpha}\bigr)^{V[G]}$, so $(C^{2b})^V = (C^{2b})^{V[G]}$. Hence the inner model $C^{2b}$ is invariant under set forcings.\label{mainthm:item1}

\item For all partial orders $\mathbb{P} \in C^{2b}$ and all $\mathbb{P}$-generic filters $G$ over $V$, we have $(H_{\omega_1}, \in )^{C^{2b}[G]} \prec (H_{\omega_1} , \in)^{V[G]}$. Hence $\mathsf{PD}$ holds in any set generic extension of $C^{2b}$ as well.\label{mainthm:item2}

\item Let $M$ be an inner model of $\mathsf{ZF}$ such that for all partial orders $\mathbb{P} \in M$ and all $\mathbb{P}$-generic filters $G$ over $V$, we have $(H_{\omega_1}, \in )^{M[G]} \prec (H_{\omega_1} , \in)^{V[G]}$. Then $(C^{2b})^M = C^{2b}$. In particular, $(C^{2b})^{C^{2b}} = C^{2b}$ and hence the model $C^{2b}$ satisfies the Axiom of Choice. Also, $C^{2b}$ is the smallest inner model $M$ of $\mathsf{ZF}$ such that for all partial orders $\mathbb{P} \in M$ and all $\mathbb{P}$-generic filters $G$ over $V$, we have $(H_{\omega_1}, \in )^{M[G]} \prec (H_{\omega_1} , \in)^{V[G]}$.\label{mainthm:item3}

\item $C^{2b}$ is the least inner model of $\mathsf{ZFC}$ closed under $\M_n^{\#}$ operators for all $n< \omega$. $C^{2b}$ is also a premouse with no measurable cardinals, hence $C^{2b}$ satisfies $\mathsf{GCH}$, $\diamondsuit_{\kappa}$ for all uncountable regular cardinals $\kappa$, and $\square_{\kappa}$ for all uncountable cardinals $\kappa$.\label{mainthm:item4} 
\end{enumerate}
\end{Thm}

\begin{proof}
\begin{description}
\item[\ref{mainthm:item1}]
The following claim is the key point:
\begin{Claim}\label{mainthm:claim}
Let $A$ be a transitive set and $\gamma = |A|$. Suppose that $\gamma^{<\omega} \subseteq A$. Then for all set generic extensions $V[G]$ of $V$, we have $\bigl( \text{Def}_{\mathrm{L}^{2b}} (A)\bigr)^V = \bigl( \text{Def}_{\mathrm{L}^{2b}} (A)\bigr)^{V[G]}$.
\end{Claim}

\begin{proof}[Proof of Claim~\ref{mainthm:claim}]
Let $V[G]$ be any set generic extension of $V$. We show the both inclusions $\bigl( \text{Def}_{\mathrm{L}^{2b}} (A)\bigr)^V \subseteq \bigl( \text{Def}_{\mathrm{L}^{2b}} (A)\bigr)^{V[G]}$ and $\bigl( \text{Def}_{\mathrm{L}^{2b}} (A)\bigr)^V \supseteq \bigl( \text{Def}_{\mathrm{L}^{2b}} (A)\bigr)^{V[G]}$.

\begin{description}
\item[$\subseteq $] 
Let $X$ be any set in $\bigl( \text{Def}_{\mathrm{L}^{2b}} (A)\bigr)^V$. We see that $X$ is also in $\bigl( \text{Def}_{\mathrm{L}^{2b}} (A)\bigr)^{V[G]}$.

Since $X$ is in $\bigl( \text{Def}_{\mathrm{L}^{2b}} (A)\bigr)^V$, there are a second-order formula $\phi$ and $\vec{a} \in A^{<\omega}$ such that the pair $(\phi , \vec{a})$ is suitable to $(A,\in)$ and $X = \{ x \in A \mid \text{$\phi [x,\vec{a}]$ is Boolean-valid with respect to $(A,\in)$} \}$.  
By Remark~\ref{rem:boolean-validity}, if $\phi [x,\vec{a}]$ is Boolean-valid with respect to $(A,\in)$ in $V$, then so is in $V[G]$. Hence the pair $(\phi , \vec{a})$ is suitable to $(A,\in)$ in $V[G]$ as well and $X = \{ x \in A \mid \text{$\phi [x,\vec{a}]$ is Boolean-valid with respect to $(A,\in)$} \}$ in $V[G]$. Therefore, $X$ is in  $\bigl( \text{Def}_{\mathrm{L}^{2b}} (A)\bigr)^{V[G]}$, as desired.

 \item[$\supseteq$] 
 Let $X$ be any set in $\bigl( \text{Def}_{\mathrm{L}^{2b}} (A)\bigr)^{V[G]}$. We see that $X$ is also in $\bigl( \text{Def}_{\mathrm{L}^{2b}} (A)\bigr)^V$.

Since $X$ is in $\bigl( \text{Def}_{\mathrm{L}^{2b}} (A)\bigr)^{V[G]}$, there are a second-order formula $\phi$ and $\vec{a} \in A^{<\omega}$ such that in $V[G]$, the pair $(\phi , \vec{a})$ is suitable to $(A,\in)$ and $X = \{ x \in A \mid \text{$\phi [x,\vec{a}]$ is Boolean-valid with respect to $(A,\in)$} \}$.  

In general, the pair $(\phi , \vec{a})$ is {\it not} suitable to $(A,\in)$ in $V$.\footnote{For example, $\phi$ could state that \lq\lq $A$ is countable", and $A$ could be uncountable in $V$ while it is countable in $V[G]$.} 
The idea is that the statement \lq\lq $\phi [x ,\vec{a}]$ holds after collapsing $A$ to be countable" is suitable to $(A, \in)$ in $V$: Recall that $\gamma = |A|$ and $\gamma^{<\omega} \subseteq A$. 
\begin{Subclaim}\label{mainthm:subclaim}
There is a second-order formula $\psi$ such that for all $x \in A$, $N^V_A \vDash \psi [x , \vec{a}]$ if and only if for all $\text{Col}(\omega , \gamma)$-generic filters $H$ over $V$, $N^{V[H]}_A \vDash \phi [x , \vec{a}]$. %in $V[H]$, for all second-order structures $P$ with first-order universe $A$ and $E^P = \{ (y,z) \in A^2 \mid y \in z \}$, $P \vDash \phi [x , \vec{a}]$. % for all full second-order structures $N$ in $V$ with first-order universe $A$ and $E^N = \{ (y, z ) \in A^2 \mid y \in z\}$, for all $x \in A$, $N \vDash \psi [x , \vec{a}]$ if and only if for all $\text{Col}(\omega , \gamma)$-generic filters $H$ over $V$, in $V[H]$, for all second-order structures $P$ with first-order universe $A$ and $E^P = \{ (y,z) \in A^2 \mid y \in z \}$, $P \vDash \phi [x , \vec{a}]$. 
\end{Subclaim}
\begin{proof}[Proof of Subclaim~\ref{mainthm:subclaim}]

Let $\mathbb{P} = \text{Col}(\omega , \gamma)$ and $\mathbb{B} = \wp (\mathbb{P}) / \! \sim $ be the completion of $\mathbb{P}$. 
Since $\gamma^{< \omega} \subseteq A$, $\mathbb{P}$ is a subset of $A$ and $\mathbb{B}$ is a complete Boolean algebra definable over $N^V_A$. %any full second-order structure $N$ of the form $N = \bigl(A , \in , \wp (A) , \ldots \bigr)$. 

Let $M$ be a Boolean-valued $\mathcal{L}$-structure (where $\mathcal{L} = \{ E \}$) of the form $M = (A , \mathbb{B} , E^M)$ where for all $y,z \in A$, $E^M (y,z) = 1$ if $y \in z$ and $E^M (y,z) = 0$ if $y\notin z$. 
Since $\mathbb{B}$ is definable over $N^V_A$, %any full second-order structure $N$ of the form $N = (A, \in , \wp (A) , \ldots )$, 
by the definition of $\llbracket\phi [\vec{a}, \vec{f}]\rrbracket^M$ %$\Qp{\phi [\vec{a}, \vec{f}]}^M$ 
in Definition~\ref{Boolean-valued interpretation}, there is a second-order formula $\psi$ such that for all $x \in A$, 
$N^V_A\vDash \psi [x , \vec{a}]$ if and only if $\Qp{\phi [x, \vec{a}]}^M = 1_{\mathbb{B}}$. By Lemma~\ref{boolean semantics}, this $\psi$ is the desired formula. 

%and all full second-order structures $N$ of the form $N = (A, \in , \wp (A), \ldots)$, $N\vDash \psi [x , \vec{a}]$ if and only if $\Qp{\phi [x, \vec{a}]}^M = 1_{\mathbb{B}}$. By Lemma~\ref{boolean semantics}, this $\psi$ is the desired formula. 

This completes the proof of Subclaim~\ref{mainthm:subclaim}.
\end{proof}

Let $\psi$ be as in Subclaim~\ref{mainthm:subclaim}. %We claim that the pair $(\psi , \vec{a})$ is suitable to $(A,\in)$ in $V[G]$. Since in $V[G]$, the pair $(\phi , \vec{a})$ is suitable to $(A, \in)$, by Lemma~\ref{boolean semantics} in $V[G]$, for all $x \in A$, either for all set generic extensions $V[G][I]$ of $V[G]$, in $V[G][I]$, for all full second-order structures $P$ with first-order universe $A$ and $E^P = \{ (y, z) \in A^2 \mid y \in z\}$, $P \vDash \phi [x,\vec{a}]$, or for all set generic extensions $V[G][I]$ of $V[G]$, in $V[G][I]$, for all full second-order structures $P$ with first-order universe $A$ and $E^P = \{ (y, z) \in A^2 \mid y \in z\}$, $P \vDash \neg \phi [x,\vec{a}]$. 
%In particular, for all $x \in A$, either for all generic extensions $V[G][I]$ of $V[G]$, for all full second-order structures $P$ with first-order universe $A$ and $E^P = \{ (y, z) \in A^2 \mid y \in z\}$, $P \vDash \psi [x,\vec{a}]$, or all set generic extensions $V[G][I]$ of $V[G]$, in $V[G][I]$, for all full second-order structures $P$ with first-order universe $A$ and $E^P = \{ (y, z) \in A^2 \mid y \in z\}$, $P \vDash \neg \psi [x,\vec{a}]$. 
%Hence by Lemma~\ref{boolean semantics} in $V[G]$ again, the pair $(\psi , \vec{a})$ is suitable to $(A,\in)$ in $V[G]$.
We claim that the pair $(\psi , \vec{a})$ is suitable to $(A, \in)$ in $V$.  
By Remark~\ref{rem:boolean-validity} in $V$ and Subclaim~\ref{mainthm:subclaim}, it is enough to show that for all $x \in A$, either for all set generic extensions $V[I]$ of $V$ and all $\text{Col}(\omega , \gamma)$-generic filters $H$ over $V[I]$, we have $N^{V[I][H]}_A \vDash \phi [x, \vec{a}]$, or for all set generic extensions $V[I]$ of $V$ and all $\text{Col}(\omega , \gamma)$-generic filters $H$ over $V[I]$, we have $N^{V[I][H]}_A \vDash \neg \phi [x , \vec{a}]$. %for all full second-order structures $P$ with first-order universe $A$ and $E^P = \{ (y,z) \in A^2 \mid y \in z\}$, we have $P \vDash \phi [x, \vec{a}]$, or for all set generic extensions $V[I]$ of $V$ and all $\text{Col}(\omega , \gamma)$-generic filters $H$ over $V[I]$, in $V[I][H]$, for all full second-order structures $P$ with first-order universe $A$ and $E^P = \{ (y,z) \in A^2 \mid y \in z\}$, we have $P \vDash \neg \phi [x, \vec{a}]$.  

Fix an $x \in A$. Since the pair $(\phi , \vec{a})$ is suitable to $(A, \in)$ in $V[G]$, $V[G]$ thinks that either $\phi [x,\vec{a}]$ is Boolean-valid with respect to $(A,\in)$, or $\neg \phi [x,\vec{a}]$ is Boolean-valid with respect to $(A,\in)$. We may assume that $\phi [x,\vec{a}]$ is Boolean-valid with respect to $(A, \in)$ in $V[G]$ (the other case is similarly argued). 
We will show that for all $I$ and $H$ as in the last paragraph, $N^{V[I][H]}_A \vDash \phi [x,\vec{a}]$. 

Fix $I$ and $H$. We may assume that $I$ and $H$ are sufficiently generic so that $V[G][I]$ is a set generic extension of $V[G]$ and $V[G][I][H]$ is a set generic extension of $V[G][I]$ with $\text{Col}(\omega , \gamma)$. 
Since $\phi [x,\vec{a}]$ is Boolean-valid with respect to $(A, \in)$ in $V[G]$, in $V[G][I][H]$, $N^{V[G][I][H]}_A \vDash \phi [x, \vec{a}]$. % for all full second-order structures $Q$ with first-order universe $A$ and $E^Q = \{ (y,z) \in A^2 \mid y \in z \}$, we have $Q \vDash \phi [x, \vec{a}]$. 

Now since $|A| = \gamma$ in $V$ and $H$ is $\text{Col}(\omega , \gamma)$-generic over $V$, $A$ is countable and transitive in $V[H]$. So $N^{V[I][H]}_A$ is first-order definable over the structure $(H_{\omega_1} , \in)^{V[I][H]}$ while $N^{V[G][I][H]}_A$ is first-order definable over $(H_{\omega_1} , \in)^{V[G][I][H]}$ by the same formula. Since $\mathsf{PD}$ holds in any set generic extension of $V[I][H]$, by Theorem~\ref{thm:projective-absoluteness}~\ref{projective-absoluteness-item3}, $(H_{\omega_1} , \in)^{V[I][H]}$ is an elementary substructure of $(H_{\omega_1}, \in)^{V[I][H][G]} = (H_{\omega_1}, \in)^{V[G][I][H]}$. Since $N^{V[G][I][H]}_A \vDash \phi [x, \vec{a}]$, we have $N^{V[I][H]}_A \vDash \phi [x, \vec{a}]$, as desired. 
%Since $P = (A ,\in , \wp (A) , \ldots )$ is in $V[I][H]$, and $Q = (A, \in , \wp (A) , \ldots )$ is in $V[G][I][H]$, $P$ is first-order definable over $(H_{\omega_1} , \in)^{V[I][H]}$ while $Q$ is first-order definable over $(H_{\omega_1} , \in)^{V[G][I][H]}$ by the same formula. Since $\mathsf{PD}$ holds in any set generic extension of $V$, by Theorem~\ref{thm:projective-absoluteness}, $(H_{\omega_1} , \in)^{V[I][H]}$ is an elementary substructure of $(H_{\omega_1}, \in)^{V[I][H][G]} = (H_{\omega_1}, \in)^{V[G][I][H]}$. Since $Q \vDash \phi [x, \vec{a}]$, we have $P \vDash \phi [x, \vec{a}]$, as desired. 

The arguments in the last three paragraphs show that for all $x \in A$, if $\phi [x,\vec{a}]$ is Boolean-valid with respect to $(A, \in)$ in $V[G]$, then $\psi [x,\vec{a}]$ is Boolean-valid with respect to $(A, \in )$ in $V$, and if $\neg \phi [x ,\vec{a}]$ is Boolean-valid with respect to $(A, \in)$ in $V[G]$, then $\neg \psi [x,\vec{a}]$ is Boolean-valid with respect to $(A, \in)$ in $V$. Since $X = \{ x \in A \mid \text{$\phi [x,\vec{a}]$ is Boolean-valid with} $ $\text{respect to $(A,\in)$ in $V[G]$}\}$, we have $ X =  \{ x \in A \mid \text{$\psi [x,\vec{a}]$ is Boolean-}$ $\text{valid with respect to $(A,\in)$ in $V$}\}$. Hence $X$ is in $\text{Def}_{\mathrm{L}^2b} (A)$, as desired. 
%We now finish the proof of Sublaim~\ref{mainthm:subclaim} by verifying that $X = \{ x \in A \mid \psi [x, \vec{a}]
%This completes the proof of Subclaim~\ref{mainthm:subclaim}.
\end{description}

%To be filled in.

This completes the proof of Claim~\ref{mainthm:claim}.
\end{proof}

We now finish the proof of \ref{mainthm:item1}. 
Let $V[G]$ be any set generic extension of $V$. Notice that it is easy to verify by induction on $\alpha$ that for all ordinals $\alpha$, $\L^{2b}_{\alpha}$ is transitive. 
We show that for all ordinals $\alpha$, $\bigl(\mathrm{L}^{2b}_{\alpha}\bigr)^V = \bigl(\mathrm{L}^{2b}_{\alpha}\bigr)^{V[G]}$ by induction on $\alpha$. 

When $\alpha < \omega$, it is easy to see that $\bigl(\mathrm{L}^{2b}_{\alpha}\bigr)^V = \mathrm{L}_{\alpha} = \bigl(\mathrm{L}^{2b}_{\alpha}\bigr)^{V[G]}$.
 
When $\alpha$ is a limit ordinal, it is clear that $\bigl(\mathrm{L}^{2b}_{\alpha}\bigr)^V = \bigl(\mathrm{L}^{2b}_{\alpha}\bigr)^{V[G]}$ by the definition of $\mathrm{L}^{2b}_{\alpha}$ and induction hypothesis.

Suppose that $\alpha = \beta +1$ for some $\beta \ge \omega$. 
By induction hypothesis, $\bigl(\mathrm{L}^{2b}_{\beta}\bigr)^V = \bigl(\mathrm{L}^{2b}_{\beta}\bigr)^{V[G]}$. 
Let $A = \bigl(\mathrm{L}^{2b}_{\beta}\bigr)^V = \bigl(\mathrm{L}^{2b}_{\beta}\bigr)^{V[G]}$ and set $\gamma = |A|$. Then $\gamma$ is a limit ordinal with $\gamma \le \beta < \alpha$. 
By the construction of $\mathrm{L}^{2b}_{\gamma}$, it is easy to see that $\gamma^{<\omega} \subseteq \bigl(\mathrm{L}^{2b}_{\gamma}\bigr)^V \subseteq \bigl(\mathrm{L}^{2b}_{\beta}\bigr)^V = A$.  
Also, $A = (\L^{2b}_{\beta})^V$ is transitive. Hence, by Claim~\ref{mainthm:claim}, $\bigl(\mathrm{L}^{2b}_{\alpha}\bigr)^V = \bigl(\mathrm{L}^{2b}_{\beta +1}\bigr)^V = \bigl( \text{Def}_{\mathrm{L}^{2b}} (A)\bigr)^V = \bigl( \text{Def}_{\mathrm{L}^{2b}} (A)\bigr)^{V[G]} = \bigl(\mathrm{L}^{2b}_{\beta +1}\bigr)^{V[G]} = \bigl(\mathrm{L}^{2b}_{\alpha}\bigr)^{V[G]}$, as desired.
 
 This completes the proof of \ref{mainthm:item1}.
 
\item[\ref{mainthm:item2}]

Let $\mathbb{P}$ be a partial order in $C^{2b}$ and $G$ be $\mathbb{P}$-generic over $V$. We show that $(H_{\omega_1}, \in)^{C^{2b}[G]}$ is an elementary substructure of $(H_{\omega_1}, \in)^{V[G]}$.

Let us fix a surjection $\pi \colon \wp (\omega) \to H_{\omega_1}$ which is $\Delta_1$ over $(H_{\omega_1}, \in)$ such that $\pi^{V[G]}$ is also surjective onto $H_{\omega_1}^{V[G]}$. 
Fix $n < \omega$. Since the $\Sigma_n$-satisfaction relation over $(H_{\omega_1}, \in)$ is projective in the codes via $\pi$, by Theorem~\ref{thm:projective-uniformization}, there is a projective function $f_n$ such that $f_n$ induces a $\Sigma_n$-Skolem function for $(H_{\omega_1}, \in)$ via $\pi$, i.e., the map sending $\pi (x)$ to $\pi \bigl( f_n (x) \bigr)$ is a $\Sigma_n$-Skolem function for $(H_{\omega_1}, \in )$, and for all reals $x$ and $y$, if $\pi (x) = \pi (y)$, then $\pi \bigl( f_n (x)\bigr) = \pi \bigl( f_n (y) \bigr)$. Furthermore, $f_n^{V[G]}$ induces a $\Sigma_n$-Skolem function for $(H_{\omega_1}, \in)^{V[G]}$ in the same manner and for all reals $x$ in $V$, $f_n^V (x) = f_n^{V[G]} (x)$. 

To see that $(H_{\omega_1}, \in)^{C^{2b}[G]}$ is an elementary substructure of $(H_{\omega_1}, \in)^{V[G]}$, it is enough to show that $\wp (\omega)^{C^{2b}[G]}$ is closed under $f_n^{V[G]}$ for all $n < \omega$. 
%Since every element of $H_{\omega_1}$ is simply coded by a real (or a subset of $\omega$), it is enough to verify that the second-order structure $\bigl(\omega , \wp (\omega) , \in \bigr)^V$ is elementary sub
%By Theorem~\ref{thm:projective-uniformization}  in $V[G]$, it is enough to show that for all $n \in \omega$, the structure $(H_{\omega_1}, \in)^{C^{2b}[G]}$ is closed under the $\Sigma_n$-Skolem function $f_n^{V[G]}$ for $(H_{\omega_1}, \in )^{V[G]}$. 
It suffices to see the following: 
\begin{Claim}\label{claim:L2bSkolem}
For any $\mathbb{P}$-name $\sigma$ for a real in $C^{2b}$, there is a $\mathbb{P}$-name $\tau$ for a real in $C^{2b}$ such that in $V$, we have $\Vdash_{\mathbb{P}} \lq\lq \tau = f_n (\sigma) $''.
\end{Claim}

\begin{proof}[Proof of Claim~\ref{claim:L2bSkolem}]

Take an ordinal $\alpha$ with $\mathbb{P} , \sigma \in \mathrm{L}^{2b}_{\alpha}$ and $\alpha \ge \omega$. 
Let $\mathbb{B} = \wp (\mathbb{P}) / \! \sim$ be the completion of $\mathbb{P}$ in $V$. Then since $\mathbb{P} \subseteq \mathrm{L}^{2b}_{\alpha}$, the complete Boolean algebra $\mathbb{B}$ is definable over $N^V_{\mathrm{L}^{2b}_{\alpha}}$. For simplicity, we write $N^V$ for $N^V_{\mathrm{L}^{2b}_{\alpha}}$ and $N^{V[G]}$ for $N^{V[G]}_{\mathrm{L}^{2b}_{\alpha}}$.%the full second-order structure $N = \bigl(\mathrm{L}^{2b}_{\alpha} , \in \wp (\mathrm{L}^{2b}_{\alpha})\bigr)$. 

Let $\tau$ be the following:
\begin{align*}
\tau = \{ (\check{m} , p) \mid m \in \omega , p \in \mathbb{P} , \text{ and } p \Vdash \lq\lq \check{m} \in f_n (\sigma)\text{''} \} %p \le \Qp{\check{m} \in f_n (\sigma)}_{\mathbb{B}}\}.
\end{align*}
Then it is easy to verify that $\tau$ is a $\mathbb{P}$-name for a real such that $\Vdash_{\mathbb{P}} \lq\lq \tau = f_n (\sigma)$'' holds in $V$.

We show that $\tau$ is in $\mathrm{L}^{2b}_{\alpha +1}$ (so in $C^{2b}$ as well) by verifying that for all $m \in \omega$ and $p \in \mathbb{P}$, the statement \lq\lq $(\check{m}, p) \in \tau$'' is suitable to $(\mathrm{L}^{2b}_{\alpha} , \in)$. 

Let $M$ be a Boolean-valued $\mathcal{L}$-structure (where $\mathcal{L} = \{ E \}$) of the form $M = (\mathrm{L}^{2b}_{\alpha} , \mathbb{B} , E^M)$ where for all $y,z \in \mathrm{L}^{2b}_{\alpha}$, $E^M (y,z) = 1$ if $y \in z$ and $E^M (y,z) = 0$ if $y\notin z$. 
Since $f_n^{V[G]}$ is projective in $V[G]$, the statement \lq\lq $m \in f_n^{V[G]} (\sigma^G)$'' is definable over the structure $N^{V[G]}$ with parameter $\sigma^G$. %$\bigl(\mathrm{L}^{2b}_{\alpha} , \wp (\mathrm{L}^{2b}_{\alpha}) , \in \bigr)^{V[G]}$ with parameter $\sigma^G$. 
So letting $g_{\sigma} \colon \omega \to \mathbb{B}$ be such that for all $m < \omega$, $g(m) = \bigvee \{ p \in \mathbb{B} \mid (\check{n}, p) \in \sigma\}$, there is a second-order formula $\phi$ such that for all $m< \omega$, we have $\Qp{\phi [m, g_{\sigma} ]}^M = \Qp{\check{m} \in f_n (\sigma)}_{\mathbb{B}}$.

Since $\mathbb{B}$ is definable over $N^V$, %the full second-order structure $N = \bigl(\mathrm{L}^{2b}_{\alpha} , \in \wp (\mathrm{L}^{2b}_{\alpha})\bigr)$, 
by the definition of $\llbracket\phi [\vec{a}, \vec{f}]\rrbracket^M$ %$\Qp{\phi [\vec{a}, \vec{f}]}^M$ 
in Definition~\ref{Boolean-valued interpretation}, there is a second-order formula $\psi$ such that for all $m < \omega$ and $p \in \mathbb{P}$, $N \vDash \psi [m, p, \sigma]$ if and only if $p \le \Qp{\phi [m, g_{\sigma} ]}^M$. 
Hence the statement $\lq\lq p \Vdash \check{m} \in f_n (\sigma)$'' %$\lq\lq p \le \Qp{\check{m} \in f_n (\sigma)}_{\mathbb{B}}$'' 
can be expressed by the second-order formula $\psi$ over $N$ and $\tau = \{ (\check{m}, p) \mid N \vDash \psi [m,p, \sigma] \}$. 
%such that for all $x \in A$ and all full second-order structures $N$ of the form $N = (A, \in , \wp (A), \ldots)$, $N\vDash \psi [x , \vec{a}]$ if and only if $\Qp{\phi [x, \vec{a}]}^M = 1_{\mathbb{B}}$. By Lemma~\ref{boolean semantics}, this $\psi$ is the desired formula. 

We now verify that the statement \lq\lq $(\check{m}, p) \in \tau$'' is suitable to $(\mathrm{L}^{2b}_{\alpha} , \in)$ by showing that the pair $(\psi , \sigma)$ is suitable to $(\mathrm{L}^{2b}_{\alpha} ,\in )$. 

Fix $m \in \omega$ and $p \in \mathbb{P}$. We show that either $\psi [m, p, \sigma]$ is Boolean-valid with respect to $(\mathrm{L}^{2b}_{\alpha} , \in)$ or $\neg \psi [m, p, \sigma]$ is Boolean-valid with respect to $(\mathrm{L}^{2b}_{\alpha}, \in~\!)$. 
We may assume $p \Vdash \lq\lq \check{m} \in f_n (\sigma)$'', i.e., $N \vDash \psi [m,p, \sigma]$ (the case $p \nVdash \lq\lq \check{m} \in f_n (\sigma)$'' %$p \nleq \Qp{\check{m} \in f_n (\sigma)}_{\mathbb{B}}$
is similarly argued). We argue that $\phi [m, p, \sigma]$ is Boolean-valid with respect to $(\mathrm{L}^{2b}_{\alpha} , \in)$. 

Since %$p \le \Qp{\check{m} \in f_n (\sigma)}_{\mathbb{B}}$ is equivalent to 
$p \Vdash \lq\lq \check{m} \in f_n (\sigma)$'', by Remark~\ref{rem:boolean-validity}, it is enough to see that for all set generic extensions $V[H]$ of $V$, $p \Vdash \lq\lq \check{m} \in f_n (\sigma)$'' holds in $V[H]$, i.e., for all $\mathbb{P}$-generic filters $G'$ over $V[H]$ with $p \in G'$, we have $m \in f_n^{V[H][G']} (\sigma^{G'})$.  

Take any $\mathbb{P}$-generic filter $G'$ over $V[H]$. Since %$p \le \Qp{\check{m} \in f_n (\sigma)}_{\mathbb{B}}$ in $V$, we have 
$p \Vdash \check{m} \in f_n (\sigma)$'' in $V$, we have $m \in f_n^{V[G']} (\sigma^{G'})$. Since $f_n^{V[G']}$ is projective in the codes in $V[G']$ and $\mathsf{PD}$ holds in any set generic extension of $V[G']$, by Theorem~\ref{thm:projective-absoluteness}~\ref{projective-absoluteness-item3}, $f_n^{V[G]}$ is invariant under set forcings, so $f_n^{V[G']} (\sigma^{G'}) = f_n^{V[G'][H]} (\sigma^{G'})$, which is equal to $f_n^{V[H][G']} (\sigma^{G'})$. Hence, $m \in f_n^{V[H][G']} (\sigma^{G'})$ holds, as desired.

Therefore, the statement \lq\lq $(\check{m}, p) \in \tau$'' is suitable to $(\mathrm{L}^{2b}_{\alpha} , \in)$ and $\tau \in \L^{2b}_{\alpha + 1}$. So $\tau \in C^{2b}$, as desired.
%To be filled in.

This completes the proof of Claim~\ref{claim:L2bSkolem}.
\end{proof}
This completes the proof of \ref{mainthm:item2}.

\item[\ref{mainthm:item3}]
Let $M$ be an inner model of $\mathsf{ZF}$ such that for all partial orders $\mathbb{P} \in M$ and all $\mathbb{P}$-generic filters $G$ over $V$, we have $(H_{\omega_1}, \in )^{M[G]} \prec (H_{\omega_1} , \in)^{V[G]}$. We will show that $(C^{2b})^M = C^{2b}$.

%By the definition of $C^{2b}$ and $(\mathrm{L}^{2b}_{\alpha} \mid \alpha \in \text{Ord})$, 
As in the proof of \ref{mainthm:item1}, it is enough to verify the following: Let $A \in M$ be transitive and $\gamma = | A|$ with $\gamma^{<\omega} \subseteq A$. Then $\bigl(\text{Def}_{\mathrm{L}^{2b}} (A)\bigr)^M = \bigl(\text{Def}_{\mathrm{L}^{2b}} (A)\bigr)^V$. 

 Let $A \in M$ be transitive and $\gamma = | A|$ with $\gamma^{<\omega} \subseteq A$.  We will show that $\bigl(\text{Def}_{\mathrm{L}^{2b}} (A)\bigr)^M = \bigl(\text{Def}_{\mathrm{L}^{2b}} (A)\bigr)^V$. 

Let $n < \omega$ and $F_n \colon H_{\omega_1} \to H_{\omega_1}$ be as follows: $F_n(a) = \{ X \subseteq a \mid $ $X$ is $\Sigma^1_n$-definable over the full second-order structure $N_a^V$ with parameters in $a \}$. 
Notice that $F_n(a)$ is in $H_{\omega_1}$ because $F_n(a)$ consists of countable subsets of $a$ and $a$ is in $H_{\omega_1}$. 
It is easy to see that $F_n$ is projective in the codes via $\pi$, where $\pi \colon \wp (\omega ) \to H_{\omega_1}$ is a surjection which is $\Delta_1$ over $(H_{\omega_1}, \in)$ as in the proof of Claim~\ref{claim:L2bSkolem}. 
Hence for all partial orders $\mathbb{P}$ in $M$ and all $\mathbb{P}$-generic filters $G$ over $V$, the function $F_n$ is absolute between $M[G]$ and $V[G]$, i.e., for all $a \in H_{\omega_1}^{M[G]}$, $F_n^{M[G]} (a) = F_n^{V[G]} (a)$. In particular, if $\mathbb{P} = \text{Col} (\omega , A)$, then $A$ is in $H_{\omega_1}^{M[G]}$ because $A$ is countable and transitive in $M[G]$. So $F_n^{M[G]} (A) = F_n^{V[G]} (A)$ in this case.

Here is the key point to see: If $G$ is a $\text{Col}(\omega , A)$-generic filter over $V$, then $\bigl(\text{Def}_{\mathrm{L}^{2b}} (A)\bigr)^M = \bigcup_{n < \omega} F_n^{M[G]} (A)$ and $\bigl(\text{Def}_{\mathrm{L}^{2b}} (A)\bigr)^V = \bigcup_{n < \omega} F_n^{V[G]} (A)$, which would give us $\bigl(\text{Def}_{\mathrm{L}^{2b}} (A)\bigr)^M = \bigl(\text{Def}_{\mathrm{L}^{2b}} (A)\bigr)^V$, because $F_n^{M[G]} (A) = F_n^{V[G]} (A)$ for all $n< \omega$ as in the last paragraph.

Let $G$ be $\text{Col}(\omega , A)$-generic over $V$. We only verify that $\bigl(\text{Def}_{\mathrm{L}^{2b}} (A)\bigr)^V = \bigcup_{n < \omega} F_n^{V[G]} (A)$ (the equality $\bigl(\text{Def}_{\mathrm{L}^{2b}} (A)\bigr)^M = \bigcup_{n< \omega} F_n^{M[G]} (A)$ follows from the same argument using the fact that $\mathsf{PD}$ holds in any set generic extension of $M$). 

Let $X$ be in $\bigl(\text{Def}_{\mathrm{L}^{2b}} (A)\bigr)^V$. We argue that $X$ is also in $F_n^{V[G]}(A)$ for some $n < \omega$. 
Since $X$ is in $\bigl(\text{Def}_{\mathrm{L}^{2b}} (A)\bigr)^V$, there are a second-order formula $\phi$ and $\vec{a} \in A^{<\omega}$ such that the pair $(\phi , \vec{a})$ is suitable to $(A , \in)$ and $X = \{ x \in A \mid \phi [x,\vec{a}] \text{ is Boolean-valid}$ $\text{with respect to $(A,\in)$}\}$. For some $n < \omega$, $\phi$ is a $\Sigma^1_n$-formula. 
Since $G$ is $\text{Col}(\omega , A)$-generic over $V$, by Lemma~\ref{boolean semantics}, $X = \{ x \in N^{V[G]}_A \vDash \phi [x , \vec{a}]\}$. Hence $X$ is definable over the full second-order structure $N_A^{V[G]}$ with parameters in $A$. Since $A$ is in $H_{\omega_1}^{V[G]}$ and $\phi$ is a $\Sigma^1_n$-formula, $X$ is in $F_n^{V[G]} (A)$, as desired.

Let $X$ be in $F_n^{V[G]} (A)$ for some $n < \omega$. We argue that $X$ is also in $\bigl(\text{Def}_{\mathrm{L}^{2b}} (A)\bigr)^V$.
Since $\gamma^{<\omega} \subseteq A$ where $\gamma = |A|$, by Claim~\ref{mainthm:claim}, it is enough to see that $X$ is in $\bigl(\text{Def}_{\mathrm{L}^{2b}} (A)\bigr)^{V[G]}$. 
Since $X$ is in $F_n^{V[G]} (A)$, there are a $\Sigma^1_n$-formula $\phi$ and $\vec{a} \in A^{<\omega}$ such that $X = \{ x \in A \mid N_A^{V[G]} \vDash \phi [x, \vec{a}] \}$. 
Since $\mathsf{PD}$ holds in any set generic extension of $V[G]$ and $A$ is in $H_{\omega_1}^{V[G]}$, 
by Theorem~\ref{thm:projective-absoluteness}~\ref{projective-absoluteness-item3} and Remark~\ref{rem:boolean-validity} in $V[G]$, the pair $(\phi , \vec{a})$ is suitable to $(A, \in)$ in $V[G]$ and $X = \{ x \in A \mid \phi [x, \vec{a}] \text{ is Boolean-valid with respect to $(A,\in)$}\} \text{ in } V[G]$. 
Hence $X$ is in $\bigl(\text{Def}_{\mathrm{L}^{2b}} (A)\bigr)^{V[G]}$, as desired.

This competes the proof of \ref{mainthm:item3}.

\item[\ref{mainthm:item4}]

For the proof of this item, we assume that the reader is familiar with the basics of inner model theory (see e.g., \cite{outline-IMT} and \cite{FSIT}). Also, we use some results from~\cite{CMI}. 

We first note that for all $n < \omega$ and all transitive sets $x$, $\M_n^{\#} (x)$ exists and it is fully iterable. This follows from the assumption that $\mathsf{PD}$ holds in any set generic extension of $V$ and from Theorem~\ref{thm:PDmouse}. 

We next claim that $C^{2b}$ is the least inner model $N$ of $\mathsf{ZFC}$ that is closed under $\M_n^{\#}$ operators for all $n < \omega$, i.e., for all $n < \omega$ and all transitive sets $x$ in $N$, we have $\M_n^{\#} (x) \in N$. 

We show that the model $C^{2b}$ is closed under $\M_n^{\#}$ operators for all $n < \omega$. 
Let $n < \omega$ and $x$ be a transitive set in $C^{2b}$. 
We have $\M_n^{\#} (x)$ in $V$. 
Let $G$ be any $\text{Col}(\omega , x)$-generic filter over $V$. 
Then the statement \lq\lq $y = \M_n^{\#} (x)$'' is definable over $(H_{\omega_1}, \in)^{V[G]}$. 
By \ref{mainthm:item2}, $(H_{\omega_1}, \in)^{C^{2b}[G]}$ is an elementary substructure of $(H_{\omega_1}, \in)^{V[G]}$. Hence $\M_n^{\#} (x)$ is in $C^{2b}[G]$. Since $G$ is arbitrary and $\M_n^{\#} (x)$ does not depend on the choice of $G$, by the weak homogeneity of $\text{Col}(\omega , x)$, $\M_n^{\#} (x)$ is in $C^{2b}$, as desired. 

Notice that from the above arguments, we have also verified that $C^{2b}$ thinks \lq\lq For all $n < \omega$ and all transitive sets $x$, $\M_n^{\#} (x)$ exists''. 

Let $N$ be an inner model of $\mathsf{ZFC}$ closed under $\M_n^{\#}$ operators for all $n < \omega$. We will see that $C^{2b} \subseteq N$. 
By Theorem~\ref{thm:projective-absoluteness}~\ref{projective-absoluteness-item2}, for all partial orders $\mathbb{P}$ in $N$ and all $\mathbb{P}$-generic filters $G$ over $V$, $(H_{\omega_1}, \in)^{N[G]}$ is an elementary substructure of $(H_{\omega_1}, \in)^{V[G]}$. 
By \ref{mainthm:item3}, we have $C^{2b} \subseteq N$, as desired. 

We now claim that there is no measurable cardinal in $C^{2b}$. 
Suppose there was a measurable cardinal $\kappa$ in $C^{2b}$. 
Let $\mu$ be a normal measure on $\kappa$ in $C^{2b}$ and let $N$ be the ultrapower of $C^{2b}$ via $\mu$. Then $N \subsetneq C^{2b}$ because $N$ is definable in $C^{2b}$ and $\mu \notin N$. 
Since $C^{2b}$ and $N$ are elementary equivalent and $C^{2b}$ thinks \lq\lq For all $n < \omega$ and all transitive sets $x$, $\M_n^{\#} (x)$ exists'', so does $N$. 
But then one can show that for all $ n < \omega$ and all transitive sets $x \in N$, we have $\M_n^{\#}(x) \in N$ by induction on $n$ using Lemma~\ref{lem:Mnsharp}. %This is because the statement \lq\lq $y = x^{\#}$'' is absolute between $N$ and $V$, and given a $k < \omega$ and a transitive set $x$, the iteration strategy of $\M_{k+1}^{\#} (x)$ is guided by $\mathcal{Q}$-structures whose complexity is at most the $\M_k^{\#}$ operator. 
Hence $N$ is an inner model of $\mathsf{ZFC}$ closed under $\M_n^{\#}$ operators for all $n < \omega$. So we have $C^{2b} \subseteq N$. But we have seen that $N \subsetneq C^{2b}$. Contradiction!

Lastly, we claim that $C^{2b}$ is a premouse. 
%Let $x$ be a transitive set. 
We say a premouse $\mathcal{M}$ is {\it countably iterable} if for all countable premice $\mathcal{N}$ with an elementary embedding $\pi \colon \mathcal{N} \to \mathcal{M}$, $\mathcal{N}$ is $(\omega_1+1)$-iterable. %For a transitive set $x$, let $\Lp (x) = \bigcup \{ M \mid \text{$M$ is a countably iterable premouse over $x$ such that its final projectum is $x$} \}$. 
We work in $C^{2b}$ and define the following sequence $(\Lp_{\alpha} \mid \alpha \in \text{Ord})$ in $C^{2b}$: 
$\Lp_0 = V_{\omega}$, $\Lp_{\alpha +1} = \bigcup \{ \mathcal{M} \mid \mathcal{M}$ is a sound countably iterable premouse such that the final projectum of $\mathcal{M}$ is equal to $(\text{Ord} \cap \Lp_{\alpha})\}$, $\Lp_{\gamma} = \bigcup_{\alpha < \gamma} \Lp_{\alpha}$ for a limit ordinal $\gamma$. Set $N = \bigcup_{\alpha \in \text{Ord}} \Lp_{\alpha}$.
Then by \cite[Lemma~2.1.2]{CMI}, $N$ is a premouse closed under any mouse operator $J$ which condenses well in $C^{2b}$ such that $C^{2b}$ thinks \lq\lq For all transitive sets $x$, $J(x)$ exists''. 
In particular, $N$ is closed under $\M_n^{\#}$ operators for all $n < \omega$. But $C^{2b}$ is the least inner model of $\mathsf{ZFC}$ closed under $\M_n^{\#}$ operators for all $n < \omega$. 
So $C^{2b} \subseteq N$. Since $N$ is defined in $C^{2b}$, $N \subseteq C^{2b}$. 
Therefore, $N = C^{2b}$ and $C^{2b}$ is a premouse, as desired.

Since $C^{2b}$ is a premouse with no subcompact cardinals, by \cite[Theorem~0.1]{squares}, it satisfies $\mathsf{GCH}$, $\diamondsuit_{\kappa}$ for all uncountable regular cardinals $\kappa$, and $\square_{\kappa}$ for all uncountable cardinals $\kappa$, as desired. 
%By assumption, $\mathsf{PD}$ holds in any set generic extension of $V$, so by Theorem~\cite{thm:PDmouse}, all $k \in \omega$ and all transitive sets $x$, $\text{M}^{\#}_k (x)$ exists and is fully iterable. 

This completes the proof of \ref{mainthm:item4}.
\end{description}
This completes the proof of Theorem~\ref{mainthm}.
\end{proof}

In summary of this section, we have seen the following:
\begin{itemize}
\item The inner model $C^{2b}$ constructed from Boolean-valued second-order logic is the least inner model of $\mathsf{ZFC}$ closed under $\M_n^{\#}$ operators for all $n< \omega$ assuming that Projective Determinacy ($\mathsf{PD}$) holds in any set generic extension. Under the same assumption, the model $C^{2b}$ is invariant under set forcings, and it admits a fine structure, hence it satisfies various combinatorial principles such as $\mathsf{GCH}$, $\diamondsuit_{\kappa}$ for all uncountable regular cardinals $\kappa$, and $\square_{\kappa}$ for all uncountable cardinals $\kappa$ as G\"{o}del's L does. %absolute among all transitive proper class models of $\mathsf{ZFC}$ closed under $\M_n^{\#}$ operators for all $n <\omega$ (in particular it is invariant under set forcings). 

On the other hand, by Myhill and Scott~\cite{Myhill-Scott}, the inner model constructed from full second-order logic is equal to HOD, the class of all hereditarily ordinal definable sets. It is well-known that HOD enjoys none of the properties mentioned in the last paragraph in $\mathsf{ZFC}$ even if one assumes large cardinal axioms. 
% enjoys various nice properties such as (i) generic invariance (i.e., it is absolutely defined among all set generic extensions of $V$), (ii) it is fine-structural and satisfies $\mathsf{GCH}$ and other combinatorial principles, and (iii) it is closed under $\M_n^{\#}$ operators for all $n< \omega$ and satisfies $\mathsf{PD}$ under the assumption that $\mathsf{PD}$ holds in any set generic extension of $V$ (or there are proper class many Woodin cardinals). Under the same assumption, the inner model $C^{2b}$ is characterized as the smallest inner model $N$ of $\mathsf{ZFC}$ such that for all partial order $\mathbb{P} \in N$ and all $\mathbb{P}$-generic filters $G$ over $V$, the structure $(H_{\omega_1},\in)^{N[G]}$ is an elementary substructure of  $(H_{\omega_1},\in)^{V[G]}$ (or the smallest inner model $N$ of $\mathsf{ZFC}$ closed under $\M_n^{\#}$ operators for all $n < \omega$).
%\item On the other hand, by Myhill and Scott~\cite{Myhill-Scott}, the inner model constructed from full second-order logic is equal to HOD, the class of all hereditarily definable sets. It is well-known that HOD enjoys none of the nice properties mentioned in the last item in $\mathsf{ZFC}$ even if one assumes large cardinal axioms. 
\end{itemize}

\section{Questions}\label{sec:Q}

We close this paper by raising some open questions:
\begin{Q}\label{Q1}
Is it consistent with $\mathsf{ZFC}$ that there is a weakly generically extendible cardinal while there are no generically extendible cardinals? 
If the answer is `Yes', what is the consistency strength of such a situation with $\mathsf{ZFC}$?
\end{Q}

Related to Question~\ref{Q1}, Usuba recently proved the following: Let $\kappa$ be the least weakly generically extendible cardinal. If $\kappa$ is not generically extendible, then there is a J\'{o}nsson cardinal less than $\kappa$. 

\begin{Q}(See \cite[Question~9.18]{BVSOL})
Can one characterize the L\"{o}wenheim number $\ell^{2b}$ of the Boolean-valued second-order logic in terms of large cardinals? In particular, could $\ell^{2b}$ be smaller than the least measurable cardinal?
\end{Q}

%To be filled in.

%To be filled in.

\bibliographystyle{plain}
\bibliography{myreference}

\end{document}